\documentclass[11pt,letterpaper]{article}

\usepackage[utf8]{inputenc} 
\usepackage[T1]{fontenc}    
\usepackage{hyperref}       
\usepackage{url}            
\usepackage{booktabs}       
\usepackage{amsmath,amssymb}       
\usepackage{microtype}      
\usepackage{algorithm}
\usepackage{algorithmic}
\usepackage{mathtools}
\usepackage{caption}
\usepackage{authblk}

\usepackage{natbib}
\usepackage{geometry}   
 \usepackage{amsthm}

\title{\textbf{Distributed Proximal Splitting Algorithms}\\
\textbf{with Rates and Acceleration}}

\date{Authors' final version.\\
Published in Front. Signal Process., Jan. 2022. \url{https://doi.org/10.3389/frsip.2021.776825}}
\author{Laurent Condat\thanks{Corresponding author. Contact: see https://lcondat.github.io}}
\author{Grigory Malinovsky}
\author{Peter Richt\'arik} 
\affil{King Abdullah University of Science and Technology (KAUST), Thuwal, Saudi Arabia}

\DeclareMathOperator*{\minimize}{minimize}

\newtheorem{theorem}{Theorem}
\newtheorem{remark}{Remark}
\newtheorem{lemma}{Lemma}
\newtheorem{proposition}{Proposition}
\newtheorem{assumption}{Assumption}

\begin{document}

\maketitle

\begin{abstract}
We analyze several generic proximal splitting algorithms well suited for large-scale convex nonsmooth optimization. We derive sublinear and linear convergence results with new rates on the function value suboptimality or distance to the solution, as well as new accelerated versions, using varying stepsizes. In addition, we propose distributed variants of these algorithms, which can be accelerated as well. While most existing results are ergodic, our nonergodic results significantly broaden our understanding of primal--dual optimization algorithms.\\

\end{abstract}

\section{Introduction}

We propose new algorithms for the generic convex optimization problem:
\begin{equation}
\minimize_{x\in\mathcal{X}} \; \left\{\Psi(x) \coloneqq \frac{1}{M}\sum_{m=1}^M \!\Big(F_m(x)+H_m(K_m x)\Big) + R(x) \right\},\label{eqpbdi}
 \end{equation}
 where $M\geq 1$ is typically the number of parallel computing nodes in a distributed setting; the $K_m:\mathcal{X}\rightarrow \mathcal{U}_m$ are linear operators; $\mathcal{X}$ and $\mathcal{U}_m$ are real Hilbert spaces (all spaces are supposed of finite dimension); $R$ and $H_m$ are proper, closed, convex functions with values in $\mathbb{R}\cup\{+\infty\}$, the proximity operators of which are easy to compute; and the $F_m$ are convex  $L_{F_m}$-smooth functions; that is   $\nabla F_m$ is $L_{F_m}$-Lipschitz continuous, for some $L_{F_m}>0$.
 
This template problem covers most convex optimization problems met in  signal and image processing,  operations research, control, machine learning, 
 and many other fields, and our goal is to propose new generic distributed algorithms able to deal with nonsmooth functions using their proximity operators, with acceleration in presence of strong convexity.
 
 \subsection{Contributions}
 
Our contributions are the following:
 
\begin{enumerate}
\item 
\textbf{New algorithms:} We propose the first distributed algorithms to solve \eqref{eqpbdi} in whole generality, with proved convergence to an exact solution, 
and having  the \emph{full splitting}, or decoupling, property:  $\nabla F_m$, $\mathrm{prox}_{H_m}$, $K_m$ and $K_m^*$  are applied at the $m$-th node, and  the proximity operator of $R$ is applied at the master node connected to all others. No other more complicated operation, like an inner loop or a linear system to solve, is involved.  

\item \textbf{Unified framework:} The foundation of our distributed algorithms consists in two general principles, applied in a cascade, which are new contributions in themselves and could be used in other contexts:

\begin{enumerate}
\item  We show that problem \eqref{eqpbdi} with $M=1$, i.e.\ the minimization of $F+R+H\circ K$, can be reformulated as the minimization of $\widetilde{F}+\widetilde{R}+\widetilde{H}$ in a different space, with preserved smoothness and strong convexity properties. Hence, the linear operator disappears and the Davis--Yin algorithm \citep{dav17} can be applied to this new problem. Through this lens, we recover many algorithms as particular cases of this unified framework, like the PD3O, Chambolle--Pock, Loris--Verhoeven algorithms.

\item We design a non-straightforward lifting technique, so that the problem \eqref{eqpbdi}, with any $M$, is reformulated as the minimization of $\widehat{F}+\widehat{R}+\widehat{H}\circ \widehat{K}$ in some product space. 
\end{enumerate}

\item \textbf{New convergence analysis and acceleration:} Even when $M=1$, we improve upon the state of the art in two ways:
\begin{enumerate}
\item For constant stepsizes, we recover existing algorithms, but we provide new, more precise, results about their convergence speed, see Theorems~\ref{th1} and \ref{th4}.
\item With a particular strategy of varying stepsizes, we exhibit new algorithms, which are accelerated versions of them. We prove $O(1/k^2)$ convergence rate on the last iterate, see Theorems~\ref{th2} and \ref{th3}, whereas  current results in the literature are ergodic, e.g.\ \cite{cha162}.
  \end{enumerate}
 \end{enumerate}

\subsection{Related Work}

Many estimation problems in a wide range of scientific fields  can be formulated as large-scale convex optimization problems~\citep{pal09, sra11, bac12, pol15, bub15, glo16, cha16, sta16,con173,con192}.  Proximal splitting algorithms~\citep{com10,bot14,par14,kom15,bec17,con19} are particularly well suited to solve them; 
they consist of simple, easy to compute, steps that can 
deal with the terms in the 
objective function separately. 

These algorithms are generally designed as sequential ones, for $M=1$, and then they can be extended by lifting in product space  
to parallel versions, well suited  to minimize $F+R+\sum_m H_m \circ K_m$, see for instance  \citealp[Section 8]{con19}. However, it is not straightforward to adapt lifting to the case of a finite-sum $F=\frac{1}{M} \sum_m F_m$, with each function $F_m$ handled by a different node, which is of primary importance in machine learning. This generalization is one of our contributions. 

There is a vast literature on distributed optimization  to minimize 
 $\frac{1}{M} \sum_m F_m + R$, with a focus on   strategies based on (block-)coordinate or randomized activation, as well as replacing the gradients by cheaper stochastic estimates~\citep{cev14,ric14,gor20,sal20}. 
Replacing the full gradient by a stochastic oracle in the accelerated algorithms with varying stepsizes we propose is not straightforward; we leave this direction for future research. In any case, the generalized setting, with the smooth  functions $F_m$ at the nodes supplemented or replaced by  nonsmooth functions $H_m$, possibly composed with linear operators, seems to have received little attention. We want to make up for that. Decentralized optimization over networks is an active research topic \citep{lat19,alg21}. In this paper, we focus on the centralized client--server model, with one master node connected to several client nodes, working in parallel. We leave the study of decentralized algorithms for future work.

When $M=1$ and $K=I$,  where $I$ denotes the identity, Davis and Yin \citep{dav17} proposed an efficient algorithm, along with an extensive study of its convergence rates and possible accelerations. 
But the ability to handle a nontrivial $K$ is behind the success of the Chambolle--Pock~\citep{cha11a} or Condat--V\~u algorithms~\citep{con13,vu13}: they are well suited for regularized inverse problems in imaging \citep{cha16}, for instance with the total variation and its variants \citep{con14,con17,dur16,bre10}; other examples are computer vision problems \citep{cre11}, overlapping group norms for sparse estimation in data science \citep{bac12}, and trend filtering on graphs \citep{wan16}. Another prominent case is when $H$ is an indicator function, so that the problem becomes: minimize $F(x)+R(x)$ subject to $Kx=b$. If $K$ is a gossip matrix like the minus graph Laplacian, decentralized optimization over a network can be tackled  \citep{shi15,sca17,sal21}.

When $M=1$ and $K$ is arbitrary, there exist algorithms to solve \eqref{eqpbdi} in full generality, for example, the Combettes--Pesquet \citep{com12}, Condat--V\~u \citep{con13,vu13}, PD3O  \citep{yan18} and PDDY \citep{sal20} algorithms. However, their convergence rates and possible accelerations are little understood. Our main contribution is to derive new convergence rates and accelerated versions of the PD3O and PDDY algorithms, and their particular cases, including  Chambolle--Pock~\citep{cha11a} and Loris--Verhoven \citep{lor11} algorithms. In order to do this, we show that these two algorithms can be viewed as instances of the Davis--Yin algorithm. This reformulation technique is inspired by 
the recent one of O'Connor and Vandenberghe~\citep{oco20}; 
it  makes it possible  to split the composition $H\circ K$ and to derive algorithms, which call the operators $\mathrm{prox}_H$, $K$, $K^*$ separately. This technique is fundamentally different from the one in \cite{sal20}, showing that the PD3O and PDDY algorithms are primal--dual instances of the operator version of Davis--Yin splitting to solve monotone inclusions. Notably, we can derive convergence rates with respect to the objective function and accelerations, which is not possible with the primal--dual reformulation of \cite{sal20}. On the other hand, the latter encompasses the Condat--V\~u algorithm~\citep{con13,vu13}, which is not the case of our approach. So, these are complementary interpretations.

\subsection{Organization of the paper}

In Section~\ref{secm1}, we propose new nonstationary versions (i.e. with varying stepsizes) of several algorithms for optimization problems made of three terms, and we analyze their convergence rates. The derivation details are pushed to the end of the paper in Section~\ref{secA5} for ease of reading. 
In Section~\ref{secdi}, we further propose  distributed algorithms, which can minimize the sum of an arbitrary number of terms.  Again, the derivation details are deferred to  Section~\ref{secA6}.
Numerical experiments illustrating the good match between our theoretical results and practical performance are shown in Section~\ref{secexp}.

\section{Minimization of 3 Functions with a Linear Operator}\label{secm1}

Let us focus on the problem \eqref{eqpbdi} when $M=1$:
\begin{equation}
\minimize_{x\in\mathcal{X}} \; \Psi(x)=F(x)+R(x)+H(K x),\label{eqpb}
 \end{equation}
 where $K:\mathcal{X}\rightarrow \mathcal{U}$ is a linear operator, $\mathcal{X}$ and $\mathcal{U}$ are real Hilbert spaces, $R$ and $H$ are proper, closed, convex functions, and $F$ is a convex and $L_F$-smooth function. We will see in Section~\ref{secdi} that using an adequate lifting technique, \eqref{eqpb} can be extended to \eqref{eqpbdi} and, accordingly, parallel or distributed versions of the sequential algorithms to solve \eqref{eqpb} will be derived. That is why we first study the case $M=1$. 
 For any function $G$, we denote by $\mu_G\geq 0$ some constant such that $G$ is $\mu_G$-strongly convex; that is, $G-(\mu_G/2)\|\cdot\|^2$ is convex.

The dual problem to \eqref{eqpb} is
\begin{equation}
\minimize_{u\in\mathcal{U}} \;(F+R)^*(-K^*u)+H^*(u),\label{eqpd}
 \end{equation}
 where $K^*$ is the adjoint operator of $K$ and $G^*$
 is the convex conjugate of a function $G$~\citep{bau17}; we recall the Moreau identity: $\mathrm{prox}_{\tau G}(z)=z-\tau\,\mathrm{prox}_{G^*/\tau}(z/\tau)$~\citep{bau17}. 
 We suppose that the following holds:
 \begin{assumption}
There exists $x^\star \in \mathcal{X}$ such that $0\in \nabla F(x^\star) + \partial R(x^\star) + K^*\partial H (Kx^\star)$, which implies that $x^\star$ is a solution to \eqref{eqpb}; see for instance \citealp[Proposition 4.3]{com12} for sufficient conditions on the functions for this property to hold. 
\end{assumption}

\subsection{Deriving the Nonstationary PD3O and PDDY Algorithms}

The main difficulty in  \eqref{eqpb} is the presence of the linear operator $K$. Indeed, 
if $K=I$, the Davis--Yin algorithm \citep{dav17} is well suited to minimize $F+R+H$. 
Note that there is a minor mistake in the way Algorithm 3 in \cite{dav17} is initialized. This is corrected here. Thus, the Davis--Yin algorithm is as follows:

Let $(\gamma_k)_{k\in\mathbb{N}}$ be a sequence of stepsizes. Let $x_H^0\in\mathcal{X}$ and $u^0\in\mathcal{X}$.
For $k=0,1,\ldots$ iterate 
\begin{equation}\left\lfloor\begin{array}{l}
x^{k+1} = \mathrm{prox}_{\gamma_k R}(x_H^k + \gamma_k u^k)\\
u^{k+1} = u^k+\frac{1}{\gamma_k}(x_H^k  - x^{k+1} )\\
x^{k+1}_H = \mathrm{prox}_{\gamma_{k+1} H}\big(x^{k+1} - \gamma_{k+1} u^{k+1}- \gamma_{k+1}\nabla F (x^{k+1})\big).
\end{array}\right.\label{eqdybase}\end{equation}
To make this algorithm applicable to $K\neq I$, we reformulate the problem \eqref{eqpb} as follows:

\begin{figure*}[th!]
\begin{minipage}{.482\textwidth}
	\begin{algorithm}[H]
		\caption*{\textbf{PD3O Algorithm} \ ($F+R + H\circ K$)}
		\begin{algorithmic}
			\STATE \textbf{input:} $(\gamma_k)_{k\in\mathbb{N}}$, $\eta\geq \|K\|^2$, $q^0\in\mathcal{X}$, $u^0\in\mathcal{U}$%
			\FOR{$k=0, 1, \ldots$}
			\STATE $x^{k+1} \coloneqq \mathrm{prox}_{\gamma_k R}\big(\gamma_{k}(q^k-K^*u^{k}) \big)$
			\STATE $q^{k+1}\coloneqq\frac{1}{\gamma_{k+1}}x^{k+1}  - \nabla F (x^{k+1})$
			\STATE $u^{k+1} \coloneqq \mathrm{prox}_{H^*\!/(\gamma_{k+1} \eta)} \big(u^k$\STATE $\ \ \ \ {}+\frac{1}{\eta}K(\frac{1}{\gamma_k} x^{k+1}+q^{k+1}-q^k )\big)$
			\ENDFOR
		\end{algorithmic}
	\end{algorithm}
		\begin{algorithm}[H]
		\caption*{\textbf{Davis--Yin Algorithm} \ ($F+R + H$)}
		\begin{algorithmic}
			\STATE \textbf{input:} $(\gamma_k)_{k\in\mathbb{N}}$, $s^0\in\mathcal{X}$
			\FOR{$k=0, 1, \ldots$}
			\STATE $x^{k+1} \coloneqq \mathrm{prox}_{\gamma_k R}(s^k)$
			\STATE $x_H^{k+1}\coloneqq \mathrm{prox}_{\gamma_{k+1} H}\big((1+ \frac{\gamma_{k+1}}{\gamma_k})x^{k+1}$
			\STATE $\ \ \ \ {}- \frac{\gamma_{k+1}}{\gamma_k}s^k
- \gamma_{k+1}\nabla F (x^{k+1})\big)$
			\STATE $s^{k+1} \coloneqq x^{k+1}_H + \frac{\gamma_{k+1}}{\gamma_k}(s^k - x^{k+1} )$
			\ENDFOR
		\end{algorithmic}
	\end{algorithm}
	\begin{algorithm}[H]
		\caption*{\textbf{Chambolle--Pock Algorithm I} \ ($R + H\circ K$)}
		\begin{algorithmic}
			\STATE \textbf{input:} $(\gamma_k)_{k\in\mathbb{N}}$, $\eta\geq \|K\|^2$, $x^0\in\mathcal{X}$, $u^0\in\mathcal{U}$%
			\FOR{$k=0, 1, \ldots$}
			\STATE $x^{k+1} \coloneqq  \mathrm{prox}_{\gamma_k R}\big(x^k-\gamma_{k}K^*u^{k} \big)$
			\STATE $u^{k+1}\coloneqq \mathrm{prox}_{H^*/(\gamma_{k+1} \eta)} \big(u^k+\frac{1}{\eta}K\big((\frac{1}{\gamma_{k+1}}$%
			\STATE $\ \ \ \ {}+\frac{1}{\gamma_k})x^{k+1}-\frac{1}{\gamma_{k}}x^{k} \big)\big)$
			\ENDFOR
		\end{algorithmic}
	\end{algorithm}
	\begin{algorithm}[H]
		\caption*{\textbf{Douglas--Rachford Algorithm} \ ($R + H$)}
		\begin{algorithmic}
			\STATE \textbf{input:} $(\gamma_k)_{k\in\mathbb{N}}$, $s^0\in\mathcal{X}$
			\FOR{$k=0, 1, \ldots$}
			\STATE $x^{k+1} \coloneqq \mathrm{prox}_{\gamma_k R}(s^k )$
			\STATE $x_H^{k+1}\coloneqq \mathrm{prox}_{\gamma_{k+1} H} \big((1+\frac{\gamma_{k+1}}{\gamma_k})x^{k+1}-\frac{\gamma_{k+1}}{\gamma_{k}}s^{k} \big)$
			\STATE $s^{k+1}\coloneqq x_H^{k+1}+\frac{\gamma_{k+1}}{\gamma_k}(s^{k}-x^{k+1})$
			\ENDFOR
		\end{algorithmic}
	\end{algorithm}
	\end{minipage}
	\ \ \ \ \ \ \begin{minipage}{.482\textwidth}
	\begin{algorithm}[H]
	\caption*{\textbf{PDDY Algorithm} \ ($F+R + H\circ K$)}
		\begin{algorithmic}
		\STATE \textbf{input:} $(\gamma_k)_{k\in\mathbb{N}}$, $\eta\geq \|K\|^2$, $x_R^0\in\mathcal{X}$, $u^0\in\mathcal{U}$%
			\STATE \textbf{initialize:}  $p^0\coloneqq K^* u^0$
			\FOR{$k=0, 1, \ldots$}
		\STATE $u^{k+1} \coloneqq \mathrm{prox}_{H^*/(\gamma_{k} \eta)} \big(u^k +\frac{1}{\gamma_k \eta}Kx_R^k\big)$
\STATE $p^{k+1} \coloneqq K^*u^{k+1}$
\STATE $x^{k+1}\coloneqq x_R^k  - \gamma_k (p^{k+1}-p^k)$
\STATE $x^{k+1}_R \coloneqq \mathrm{prox}_{\gamma_{k+1} R}\big(x^{k+1}-\gamma_{k+1}\nabla F (x^{k+1})$
\STATE $\ \ \ \ {}- \gamma_{k+1} p^{k+1}\big)$
			\ENDFOR
		\end{algorithmic}
	\end{algorithm}
	\begin{algorithm}[H]
	\caption*{\textbf{Loris--Verhoeven Algorithm} \ ($F+H\circ K$)}
		\begin{algorithmic}
		\STATE \textbf{input:} $(\gamma_k)_{k\in\mathbb{N}}$, $\eta\geq \|K\|^2$, $q^0\in\mathcal{X}$, $u^0\in\mathcal{U}$%
			\FOR{$k=0, 1, \ldots$}
		\STATE $x^{k+1} \coloneqq \gamma_{k}(q^k-K^*u^{k}) $
\STATE $q^{k+1} \coloneqq \frac{1}{\gamma_{k+1}}x^{k+1} 
- \nabla F (x^{k+1})$
\STATE $u^{k+1}\coloneqq  \mathrm{prox}_{H^*/(\gamma_{k+1} \eta)} \big(u^k$
\STATE $\ \ \ \ {}+\frac{1}{\eta}K(\frac{1}{\gamma_k}x^{k+1}+q^{k+1}-q^k )\big)$
			\ENDFOR
		\end{algorithmic}
	\end{algorithm}
	\begin{algorithm}[H]
		\caption*{\textbf{Chambolle--Pock Algorithm II} \ ($R + H\circ K$)}
		\begin{algorithmic}
			\STATE \textbf{input:} $(\gamma_k)_{k\in\mathbb{N}}$, $\eta\geq \|K\|^2$, $x_R^0\in\mathcal{X}$, $u^0\in\mathcal{U}$%
			\FOR{$k=0, 1, \ldots$}
			\STATE $u^{k+1} \coloneqq  \mathrm{prox}_{H^*/(\gamma_{k} \eta)} \big(u^k +\frac{1}{\gamma_k \eta}Kx_R^k\big)$
			\STATE $x^{k+1}_R \coloneqq \mathrm{prox}_{\gamma_{k+1} R}%
			\Big(x_R^k-  K^*\big((\gamma_k$%
			\STATE $\ \ \ \ {}+\gamma_{k+1}) u^{k+1}-\gamma_k u^k\big)\Big)$
			\ENDFOR
		\end{algorithmic}
	\end{algorithm}
	\begin{algorithm}[H]
		\caption*{\textbf{Forward--Backward Algorithm} \ ($F + R$)}
		\begin{algorithmic}
			\STATE \textbf{input:} $(\gamma_k)_{k\in\mathbb{N}}$,  $x_1\in\mathcal{X}$,
			\FOR{$k=1, 2,\ldots$}
			\STATE $x^{k+1} \coloneqq  \mathrm{prox}_{\gamma_{k}R} \big(x^k - \gamma_{k}\nabla F(x^k)
		\big)$
		\ENDFOR
		\end{algorithmic}
	\end{algorithm}
	\end{minipage}
	\end{figure*}

\begin{enumerate}

\item We choose a value $\eta\geq \|K\|^2$; we recommend to set $\eta=\|K\|^2$ in practice. Then there exists a real Hilbert space $\mathcal{W}$ and a linear operator $C:\mathcal{W}\rightarrow \mathcal{U}$ such that $KK^*+CC^*=\eta I$. $C$ is not unique, for instance, we can set $C=(\eta I - KK^*)^{1/2}$. We actually don't need to exhibit $C$, its existence is sufficient here and there will be no call to $C$ in the algorithms.

\item Now, the problem \eqref{eqpb} can be rewritten as:
\begin{equation}
{\minimize_{x\in\mathcal{X},w\in\mathcal{W}} \; \widetilde{F}(x,w)+\widetilde{R}(x,w)+\widetilde{H}(x,w)},\label{eqpbref}
 \end{equation}
where $\widetilde{F}:(x,w)\mapsto F(x)+\frac{\mu_F}{2}\|w\|^2$, $\widetilde{R}:(x,w)\mapsto R(x)+\imath_{0}(w)$, where $\imath_{0}:w\mapsto\{0$ if $w=0$, $+\infty$ otherwise$\}$, and $\widetilde{H}:(x,w)=H(Kx+Cw)$. Indeed, we introduce  the variable $w$, but also the constraint that $w=0$. Since $\widetilde{F}(x,0)=F(x)$, $\widetilde{R}(x,0)=R(x)$, $\widetilde{H}(x,0)=H(Kx)$, the equivalence between \eqref{eqpb} and \eqref{eqpbref} follows.
\end{enumerate}\medskip

We have $\nabla \widetilde{F}(x,w)= (\nabla F(x),\mu_F w)$, $\mathrm{prox}_{\widetilde{R}}(x,w)=(\mathrm{prox}_R(x),0)$. Most importantly, for every $\gamma>0$, we have~\citep{oco20}:
\begin{equation}
\mathrm{prox}_{\widetilde{H}^*/\gamma }(x,w)=(K^* u, C^* u ),\ \mbox{where}\  u= \mathrm{prox}_{H^*/(\gamma \eta)} \big((Kx+Cw)/\eta\big).\label{eqkkk}
\end{equation}
Note that in~\cite{oco20}, the authors use $\widetilde{F}(x,w)=F(x)$, whereas we add $\frac{\mu_F}{2}\|w\|^2$. This difference is essential, so that $\widetilde{F}$ is $L_F$-smooth and $\mu_F$-strongly convex. Also, $\widetilde{R}$ is $\mu_R$-strongly convex.

Then, we can apply the Davis--Yin algorithm \eqref{eqdybase} to solve the problem \eqref{eqpbref}. We set $F$, $R$, $H$ in \eqref{eqdybase} as $\widetilde{F}$, $\widetilde{R}$, $\widetilde{H}$, respectively. The details of the substitutions yielding the algorithms are deferred to Section~\ref{secA5} for the convenience of reading; 
most notably, whenever $CC^*$ appears, it is replaced by $\eta I-KK^*$. The obtained algorithms turns out to be a nonstationary version of the PD3O algorithm~\citep{yan18}, shown above. On the other hand, if we exchange the two functions and set $F$, $R$, $H$ in \eqref{eqdybase} as $\widetilde{F}$, $\widetilde{H}$, $\widetilde{R}$, we obtain a different algorithm. It turns out to be a nonstationary version of the PDDY algorithm proposed recently \citep{sal20}, shown above too. With constant stepsizes $\gamma_k\equiv \gamma \in (0,2/L_F)$, for both the PD3O and PDDY algorithms, $x^k$ and $u^k$ converge to some solutions $x^\star$ and $u^\star$ of \eqref{eqpb} and \eqref{eqpd}, respectively; this result was known for $\eta > \|K\|^2$~\citep{yan18,sal20} and shown for  $\eta = \|K\|^2$ for the PD3O algorithm in \cite{oco20}, but convergence with $\eta = \|K\|^2$ for the PDDY algorithm, as stated in Theorem~\ref{theoaj}, is new. \bigskip

Particular cases of the PD3O and PDDY algorithms, which are shown above, are the following:
\begin{enumerate}
\item If $K=I$ and $\eta=1$, the PD3O algorithm reverts to the Davis--Yin algorithm \eqref{eqdybase};  the PDDY algorithm too, but with $H$ and $R$ exchanged in \eqref{eqdybase}.

\item If $F=0$, the PD3O and PDDY algorithms revert to the forms I and II \citep{con19} of the Chambolle--Pock algorithm, a.k.a.\ Primal--Dual Hybrid Gradient algorithm \citep{cha11a}, respectively. 

\item If $R=0$,  the PD3O and PDDY algorithms revert to the Loris--Verhoeven algorithm~\citep{lor11}, also discovered independently as the PDFP2O~\citep{che13} and PAPC~\citep{dro15} algorithms; see also~\cite{com14,con19} for an analysis as a primal--dual forward--backward algorithm.
 
\item If $F=0$ in the Davis--Yin algorithm or $K=I$ and $\eta=1$ in the Chambolle--Pock algorithm, we obtain the Douglas--Rachford algorihm; it is equivalent to the ADMM, see the discussion in \cite{con19}.

\item If $H=0$, the PD3O and PDDY algorithms revert 
to the forward--backward algorithm, a.k.a.\ proximal gradient descent. The Loris--Verhoeven algorithm with $K=I$ and $\eta=1$, too.
\end{enumerate}

\subsection{Convergence Analysis}

We first give convergence rates for the PD3O algorithm with constant stepsizes.

\begin{theorem}[convergence rate of the PD3O algorithm]\label{th1}
In the PD3O algorithm, suppose that $\gamma_k\equiv \gamma \in (0,2/L_F)$ and $\eta\geq \|K\|^2$. Then $x^k$ and $u^k$ converge to some solutions $x^\star$ and $u^\star$  of \eqref{eqpb} and \eqref{eqpd}, respectively. In addition, suppose that  $H$ is continuous on an open ball centered at $Kx^\star$. 
Then the following hold:
\begin{equation*}
\mathrm{(i)}\quad\Psi(x^k)-\Psi(x^\star) =o(1/\sqrt{k}).
\end{equation*}
Define the weighted ergodic iterate 
$\bar{x}^{k} = \frac{2}{k(k+1)}\sum_{i=1}^{k}i x^{i}$, for every $k\geq 1$. Then 
\begin{equation*}
\mathrm{(ii)}\quad\Psi(\bar{x}^k)-\Psi(x^\star) =O(1/k). 
\end{equation*}
Furthermore, if $H$ is $L$-smooth for some $L>0$, we have a faster decay for the best iterate so far:
\begin{equation*}
\mathrm{(iii)}\quad\min_{i=1,\ldots,k}\Psi(x^i)-\Psi(x^\star) =o(1/k). 
\end{equation*}
\end{theorem}

\proof
The convergence of $x^k$ follows from \citealp[Theorem 2.1]{dav17} and the convergence of 
$u^k$ follows from the one of the variable $u_B^k=(z^k-x_A^k)/\gamma$  in the notations of \cite{dav17}.
 $\mathrm{(i)}$ follows from \citealp[Theorem 3.1]{dav17}, using the following facts; first, in this theorem, the function corresponding to $\widetilde{H}$ is supposed to be Lipschitz-continuous on a certain ball, but since the rate is asymptotic and $Kx^k\rightarrow Kx^\star$, it is sufficient to consider the property around $Kx^\star$; second, it is well known that if a convex real-valued function is continuous on a convex open set, it is Lipschitz-continuous on every compact subset of this set~\citep{unk72}; 
third, if $H$ is continuous, $\widetilde{H}$ is continuous too. 
 $\mathrm{(ii)}$ follows from \citealp[Theorem 3.2]{dav17} and $\mathrm{(iii)}$ follows from
Theorem D.5 in the preprint of \cite{dav17}.
\endproof

Theorem~\ref{th1} applies to the particular cases of the PD3O algorithm, like the Loris--Verhoeven, Chambolle--Pock, Douglas--Rachford algorithms. Our results are new even for them. 

\begin{remark}
We can note that the forward--backward algorithm $x^{k+1}=\mathrm{prox}_{\gamma R}(x^k-\gamma \nabla F(x^k))$, which is a particular case of the PD3O algorithm when $H=0$, is monotonic. So, the best iterate so far is the last iterate. Hence, Theorem~\ref{th1} $\mathrm{(iii)}$ yields $\Psi(x^k)-\Psi(x^\star) =o(1/k)$ for the forward--backward algorithm.
\end{remark}

For the PDDY algorithm, we cannot derive a similar theorem, since $\widetilde{R}$ is not continuous around $(x^\star,0)$. Still, we can establish convergence of the variables:

\begin{theorem}[convergence  of the PDDY algorithm]\label{theoaj}
In the PDDY algorithm, suppose that $\gamma_k\equiv \gamma \in (0,2/L_F)$ and $\eta\geq \|K\|^2$. Then $x^k$ and  $x_R^k$ both converge to some solution $x^\star$  of \eqref{eqpb}, and $u^k$ converges to some solution $u^\star$  of \eqref{eqpd}.
\end{theorem}
\proof
The convergence of $x^k$ and $x_R^k$ to the same solution $x^\star$ of \eqref{eqpb} follows from \citealp[Theorem 2.1]{dav17}. The convergence of the variable $u_B^k=(z^k-x_A^k)/\gamma$,  in the notations of \cite{dav17}, implies in our setting, according to \eqref{eqkkk}, that $K^*u^k$ and $C^*u^k$ both converge to some elements. But since $\eta u^k= KK^*u^k+CC^*u^k$, $u^k$ converges to some element $u^\star\in\mathcal{U}$. Finally, we have $x^\star = \mathrm{prox}_{\gamma R} (x^\star - \gamma \nabla F(x^\star)-\gamma K^* u^\star)$, so that $0\in \partial R(x^\star)+ \nabla F(x^\star) + K^* u^\star$, and $u^\star = \mathrm{prox}_{H^*/(\gamma \eta)}(u^\star + \frac{1}{\gamma\eta}Kx^\star)$, so that $Kx^\star \in (\partial H)^{-1}(u^\star)$. Hence, $u^\star$ is a solution to \eqref{eqpd}.
\endproof

We now give accelerated convergence results using varying stepsizes, when $F$ or $R$ is strongly convex; that is, $\mu_F+\mu_R>0$. In that case, we denote by $x^\star$ the unique solution to \eqref{eqpb}.

\begin{theorem}[convergence rate of the accelerated PD3O algorithm]\label{th2}
Suppose that $\mu_F+\mu_R>0$. Let $\kappa\in(0,1)$ and $\gamma_0\in (0,2(1-\kappa)/L_F)$.  Set $\gamma_1=\gamma_0$ and
\begin{equation}
\gamma_{k+1}=\frac{-\gamma_k^2\mu_F\kappa+\gamma_k\sqrt{(\gamma_k\mu_F\kappa)^2+1+2\gamma_k\mu_R}}{1+2\gamma_k\mu_R}
,\quad\mbox{for every }k\geq 1.\label{eqgammag}
\end{equation}
Suppose that $\eta\geq \|K\|^2$. Then in the PD3O algorithm, there exists $c_0>0$ (whose expression is given in Section~\ref{secA5}) such that, for every $k\geq 1$,
\begin{equation*}
\|x^{k+1}-x^\star\|^2\leq \frac{\gamma_{k+1}^2}{1-\gamma_{k+1}\mu_F\kappa}c_0
=O\big(1/k^2\big).
\end{equation*}
\end{theorem}

\proof
This result follows from \citealp[Theorem 3.3]{dav17}, stated for convenience as Lemma~\ref{lem1} in Section~\ref{secA5}.
\endproof

Note that with the stepsize rule in \eqref{eqgammag}, we have $k\,\gamma_k \rightarrow 1/(\mu_F\kappa+\mu_R)$ as $k\rightarrow+\infty$, so that $\gamma_k=O(1/k)$ and $\gamma_{k+1}/\gamma_k \rightarrow 1$. Also, when $F=0$, $L_F$ can be taken arbitrarily small, so that we can choose any $\gamma_0>0$.

Theorem~\ref{th2} is new for the PD3O and Loris--Verhoeven algorithms, but has been derived in \cite{oco20} for the Chambolle--Pock algorithm. For the forward--backward algorithm, strong convexity yields linear convergence with constant stepsizes, so this nonstationary version does not seem interesting.\bigskip

Concerning the PDDY algorithm, $\widetilde{H}$ is not necessarily strongly convex, even if $H$ is. So, we only consider the case where $F$ is strongly convex. As a consequence of Lemma~\ref{lem1}, we get:

\begin{theorem}[convergence rate of the accelerated PDDY algorithm]\label{th3}
Suppose that $\mu_F>0$. Let $\kappa\in(0,1)$ and $\gamma_0\in (0,2(1-\kappa)/L_F)$.  Set $\gamma_1=\gamma_0$ and
\begin{equation}
\gamma_{k+1}=-\gamma_k^2\mu_F\kappa+\gamma_k\sqrt{(\gamma_k\mu_F\kappa)^2+1}
,\quad\mbox{for every }k\geq 1.\label{eqgammag2}
\end{equation}
Suppose that $\eta\geq \|K\|^2$. Then in the PDDY algorithm, 
there exists $c_0>0$ (whose expression is given in Section~\ref{secA5}) such that, for every $k\geq 1$,
\begin{equation*}
\|x^{k+1}-x^\star\|^2\leq \frac{\gamma_{k+1}^2}{1-\gamma_{k+1}\mu_F\kappa}c_0
=O\big(1/k^2\big).
\end{equation*}
Moreover, if $\eta>\|K\|^2$, $\|x_R^{k}-x^\star\|^2 = O(1/k^2)$ as well.
\end{theorem}

Finally, we consider the case where, in addition to strong convexity of $F$ or $R$, $H$ is smooth; in that case, the algorithms with constant stepsizes converge linearly; that is, as a consequence of Lemma~\ref{lem2}, we have:

\begin{theorem}[linear convergence of the PD3O and PDDY algorithms]\label{th4}
Suppose that $\mu_F+\mu_R>0$ and that $H$ is $L_H$-smooth, for some $L_H>0$. Let $x^\star$ and $u^\star$ be the unique solutions to \eqref{eqpb} and \eqref{eqpd}, respectively. 
Suppose that $\gamma_k\equiv \gamma \in (0,2/L_F)$ and $\eta\geq \|K\|^2$. 
Then the PD3O algorithm converges linearly: there exists $\rho\in (0,1]$ such that, for every $k\in\mathbb{N}$,
\begin{align*}
\|x^{k+1}-x^\star\|^2\leq (1-\rho)^k &\Big(\|\gamma q^0-x^\star+\gamma \nabla F(x^\star)-\gamma K^*(u^0-u^\star)\|^2\notag\\
&\ {}+\gamma^2\eta\|u^0-u^\star\|^2-\gamma^2\|K^*(u^0-u^\star)\|^2
\Big).
\end{align*}
The PDDY algorithm converges linearly too: there exists $\rho\in (0,1]$ such that, for every $k\in\mathbb{N}$,
\begin{equation*}
\|x_R^{k+1}-x^\star\|^2\leq 4(1-\rho)^k \Big(\|x_R^0-x^\star +\gamma K^*(u^0-u^\star)\|^2
+\gamma^2\eta\|u^0-u^\star\|^2-\gamma^2\|K^*(u^0-u^\star)\|^2
\Big).
\end{equation*}
\end{theorem}

Linear convergence of the other variables in the algorithms can be derived as well, see Proposition~\ref{prop1}. 
Lower bounds for $\rho$ can be derived from Theorem D.6 in the preprint version of \cite{dav17}. We don't provide them, since they are not tight, as noticed in Remark D.2 of the same preprint.
For instance, for the PDDY or Loris--Verhoeven algorithms with $\mu_F>0$,
\begin{equation*}
\rho=\frac{\gamma \mu_F(2-\gamma L_F)}{(1+\gamma \eta L_H)^2}.
\end{equation*}
If $H=0$, by setting $L_H=0$, we get $\rho=\gamma \mu_F(2-\gamma L_F)$. But then the PDDY algorithm reverts to the  forward--backward algorithm, for which it is known 
that  $1-\rho=(1-\gamma\mu_F)^2$ whenever $\gamma\leq 2/(L_F+\mu_F)$, which corresponds to the larger value $\rho=\gamma\mu_F(2-\gamma\mu_F)$. 

We emphasize that linear convergence comes \emph{for free} with the algorithms, if the conditions are met, without any modification. That is, there is no need to know $\mu_F$, $\mu_R$, $L_H$, since the conditions on the two parameters $\gamma$ and $\eta$ do not depend on these values. For the particular case of the Chambolle--Pock algorithm, as pointed out in \cite{oco20}, this is in contrast to existing linear convergence results \citep{cha16}, derived for a modified version of the algorithm, which depends on these values.

\section{Distributed Proximal Algorithms}\label{secdi}

\begin{figure*}[t!]
\begin{minipage}{.485\textwidth}
	\begin{algorithm}[H]
		\caption*{\textbf{Distributed PD3O Algorithm}}
		\begin{algorithmic}
			\STATE \textbf{input:} $(\gamma_k)_{k\in\mathbb{N}}$, $\eta\geq \|\widehat{K}\|^2$, $(\omega_m)_{m=1}^M$,%
			\STATE \ \ \ \ $(q_m^0)_{m=1}^M\in\mathcal{X}^M$, $(u_m^0)_{m=1}^M\in\widehat{\mathcal{U}}$%
			\STATE \textbf{initialize:} $a_m^{0}\coloneqq q_m^{0}-K_m^*u_m^{0}$, $m=1...M$%
			\FOR{$k=0, 1, \ldots$}
			\STATE at master, \textbf{do}
			\STATE \ \ \ \ $x^{k+1} \coloneqq \mathrm{prox}_{\gamma_k R}\big(\frac{\gamma_{k}}{M}\sum_{m=1}^M a_m^{k} \big)$
			\STATE \ \ \ \ broadcast $x^{k+1}$ to all nodes
			\STATE at all nodes, for $m=1,\ldots,M$, \textbf{do}
			\STATE \ \ \ \ $q_m^{k+1}\coloneqq\frac{M\omega_m}{\gamma_{k+1}}x^{k+1}  - \nabla F_m (x^{k+1})$
			\STATE \ \ \ \ $u_m^{k+1} \coloneqq \mathrm{prox}_{M\omega_m H_m^*/(\gamma_{k+1} \eta)} \big(u_m^k$
			\STATE \ \ \ \ $\ \ \ \ {}+\frac{1}{\eta}K_m(\frac{M\omega_m}{\gamma_k} x^{k+1}+q_m^{k+1}-q_m^k )\big)$
			\STATE \ \ \ \ $a_m^{k+1}\coloneqq  q_m^{k+1}-K_m^*u_m^{k+1}$ 
			\STATE \ \ \ \ transmit $a_m^{k+1}$  to master
			\ENDFOR
		\end{algorithmic}
	\end{algorithm}
	\begin{algorithm}[H]
		\caption*{\textbf{Distributed Loris--Verhoeven Algorithm}}
		\begin{algorithmic}
			\STATE \textbf{input:} $(\gamma_k)_{k\in\mathbb{N}}$, $\eta\geq \|\widehat{K}\|^2$, $(\omega_m)_{m=1}^M$%
			\STATE \ \ \ \ $(q_m^0)_{m=1}^M\in\mathcal{X}^M$, $(u_m^0)_{m=1}^M\in\widehat{\mathcal{U}}$%
			\STATE \textbf{initialize:} $a_m^{0}\coloneqq q_m^{0}-K_m^*u_m^{0}$, $m=1...M$%
			\FOR{$k=0, 1, \ldots$}
			\STATE at master, \textbf{do}
			\STATE \ \ \ \ $x^{k+1} \coloneqq \frac{\gamma_{k}}{M}\sum_{m=1}^M a_m^{k}$
			\STATE \ \ \ \ broadcast $x^{k+1}$ to all nodes
			\STATE at all nodes, for $m=1,\ldots,M$, \textbf{do}
			\STATE \ \ \ \ $q_m^{k+1}\coloneqq\frac{M\omega_m}{\gamma_{k+1}}x^{k+1}  - \nabla F_m (x^{k+1})$
			\STATE \ \ \ \ $u_m^{k+1} \coloneqq \mathrm{prox}_{M\omega_m H_m^*/(\gamma_{k+1} \eta)} \big(u_m^k$
			\STATE \ \ \ \ $\ \ \ \ {}+\frac{1}{\eta}K_m(\frac{M\omega_m}{\gamma_k} x^{k+1}+q_m^{k+1}-q_m^k )\big)$
			\STATE \ \ \ \ $a_m^{k+1}\coloneqq  q_m^{k+1}-K_m^*u_m^{k+1}$ 
			\STATE \ \ \ \ transmit $a_m^{k+1}$  to master
			\ENDFOR
		\end{algorithmic}
	\end{algorithm}
	 \end{minipage}
	\ \ \ \ \ \begin{minipage}{.485\textwidth}
	\begin{algorithm}[H]
	\caption*{\textbf{Distributed PDDY Algorithm}}
		\begin{algorithmic}
		\STATE \textbf{input:} $(\gamma_k)_{k\in\mathbb{N}}$, $\eta\geq \|\widehat{K}\|^2$, $(\omega_m)_{m=1}^M$,
		\STATE \ \ \ \  $x_R^0\in\mathcal{X}$, $(u_m^0)_{m=1}^M\in\widehat{\mathcal{U}}$%
			\STATE \textbf{initialize:} $p_m^0\coloneqq K_m^* u_m^0$, $m=1,...,M$
			\FOR{$k=0, 1, \ldots$}
			\STATE at all nodes, for $m=1,\ldots,M$, \textbf{do}
		\STATE \ \ \ \ $u_m^{k+1} \coloneqq \mathrm{prox}_{M\omega_m H_m^*/(\gamma_{k} \eta)} \big(u_m^k$ 
		\STATE \ \ \ \ $\ \ \ \ {}+\frac{M\omega_m}{\gamma_k \eta}K_mx_R^k\big)$
\STATE \ \ \ \ $p_m^{k+1} \coloneqq K_m^*u_m^{k+1}$
\STATE \ \ \ \ $x_m^{k+1}\coloneqq x_R^k  - \frac{\gamma_k}{M\omega_m} (p_m^{k+1}-p_m^k)$
\STATE \ \ \ \ $a_m^k \coloneqq M\omega_m x_m^{k+1}-\gamma_{k+1}\nabla F_m (x_m^{k+1})$
\STATE \ \ \ \ $\ \ \ \ {}- \gamma_{k+1} p_m^{k+1}$
\STATE \ \ \ \ transmit $a_m^{k}$  to master
\STATE at master, \textbf{do}
\STATE \ \ \ \ $x^{k+1}_R \coloneqq \mathrm{prox}_{\gamma_{k+1} R}\big(\frac{1}{M}\sum_{m=1}^M a_m^k\big)$
\STATE \ \ \ \ broadcast $x^{k+1}_R$ to all nodes
			\ENDFOR
		\end{algorithmic}
	\end{algorithm}
	\begin{algorithm}[H]
		\caption*{\textbf{Distributed Davis--Yin Algorithm}}
		\begin{algorithmic}
			\STATE \textbf{input:} $(\gamma_k)_{k\in\mathbb{N}}$, $(s_m^0)_{m=1}^M\in\mathcal{X}^M$, $(\omega_m)_{m=1}^M$%
			\FOR{$k=0, 1, \ldots$}
			\STATE at master, \textbf{do}
			 \STATE \ \ \ \ $x^{k+1} \coloneqq \mathrm{prox}_{\gamma_k R}(\sum_{m=1}^M \omega_m s_m^k)$
			 \STATE \ \ \ \ broadcast $x^{k+1}$ to all nodes
			 \STATE at all nodes, for $m=1,\ldots,M$, \textbf{do}
			\STATE \ \ \ \ $x_{m}^{k+1}\coloneqq \mathrm{prox}_{\gamma_{k+1} H_m/(M\omega_m)} \big((1+\frac{\gamma_{k+1}}{\gamma_k})$\\
			\STATE \ \ \ \ $\ \ \ \ {}\times x^{k+1}-\frac{\gamma_{k+1}}{\gamma_{k}}s_m^{k} -\frac{\gamma_{k+1}}{M\omega_m}\nabla F_m(x^{k+1})\big)$
			\STATE \ \ \ \ $s_m^{k+1}\coloneqq x_m^{k+1}+\frac{\gamma_{k+1}}{\gamma_k}(s_m^{k}-x^{k+1})$
			\STATE \ \ \ \ transmit $s_m^{k+1}$  to master
			\ENDFOR
		\end{algorithmic}
	\end{algorithm}
	\end{minipage}
	\end{figure*}

\begin{figure*}[t!]
\begin{minipage}{.485\textwidth}
	\begin{algorithm}[H]
		\caption*{\textbf{Distributed Chambolle--Pock Algorithm}}
		\begin{algorithmic}
			\STATE \textbf{input:} $(\gamma_k)_{k\in\mathbb{N}}$, $\eta\geq \|\widehat{K}\|^2$, $(\omega_m)_{m=1}^M$
			\STATE \ \ \ \ $x_0\in\mathcal{X}$, $(u_m^0)_{m=1}^M\in\widehat{\mathcal{U}}$%
			\STATE \textbf{initialize:} $a_m^{0}\coloneqq K_m^*u_m^{0}$, $m=1,...,M$%
			\FOR{$k=0, 1, \ldots$}
			\STATE at master, \textbf{do}
			\STATE \ \ \ \ $x^{k+1} \coloneqq \mathrm{prox}_{\gamma_k R}\big(x^k-\frac{\gamma_{k}}{M}\sum_{m=1}^M a_m^{k} \big)$
			\STATE \ \ \ \ broadcast $x^{k+1}$ to all nodes
			\STATE at all nodes, for $m=1,\ldots,M$, \textbf{do}
			\STATE \ \ \ \ $u_m^{k+1} \coloneqq \mathrm{prox}_{M\omega_m H_m^*/(\gamma_{k+1} \eta)} \big(u_m^k$
			\STATE \ \ \ \ $\ \ \ \ {}+\frac{M\omega_m}{\eta}K_m\big((\frac{1}{\gamma_k}+\frac{1}{\gamma_{k+1}}) x^{k+1}-\frac{1}{\gamma_{k}}x^k \big)\big)$%
			\STATE \ \ \ \ $a_m^{k+1}\coloneqq  K_m^*u_m^{k+1}$ 
			\STATE \ \ \ \ transmit $a_m^{k+1}$  to master
			\ENDFOR
		\end{algorithmic}
	\end{algorithm}
	\begin{algorithm}[H]
		\caption*{\textbf{Distributed Douglas--Rachford Algorithm}}
		\begin{algorithmic}
			\STATE \textbf{input:} $(\gamma_k)_{k\in\mathbb{N}}$, $(\omega_m)_{m=1}^M$, $(s_m^0)_{m=1}^M\in\mathcal{X}^M$%
			\FOR{$k=0, 1, \ldots$}
			\STATE at master, \textbf{do}
			 \STATE \ \ \ \ $x^{k+1} \coloneqq \mathrm{prox}_{\gamma_k R}\big(\sum_{m=1}^M \omega_m s_m^k\big)$
			 \STATE \ \ \ \ broadcast $x^{k+1}$ to all nodes
			 \STATE at all nodes, for $m=1,\ldots,M$, \textbf{do}
			\STATE \ \ \ \ $x_{m}^{k+1}\coloneqq \mathrm{prox}_{\gamma_{k+1} H_m/(M\omega_m)} $\\ 
			\STATE \ \ \ \ $\ \ \ \ \big((1+\frac{\gamma_{k+1}}{\gamma_k})x^{k+1}-\frac{\gamma_{k+1}}{\gamma_{k}}s_m^{k} \big)$
			\STATE \ \ \ \ $s_m^{k+1}\coloneqq x_m^{k+1}+\frac{\gamma_{k+1}}{\gamma_k}(s_m^{k}-x^{k+1})$
			\STATE \ \ \ \ transmit $s_m^{k+1}$  to master
			\ENDFOR
		\end{algorithmic}
	\end{algorithm}
	\end{minipage}
\ \ \ \ \ 
\begin{minipage}{.485\textwidth}
\begin{algorithm}[H]
	\caption*{\textbf{Distributed Chambolle--Pock Alg.\ Form II}}
		\begin{algorithmic}
		\STATE \textbf{input:} $(\gamma_k)_{k\in\mathbb{N}}$, $\eta\geq \|\widehat{K}\|^2$, $(\omega_m)_{m=1}^M$,
		\STATE \ \ \ \  $x_R^0\in\mathcal{X}$, $(u_m^0)_{m=1}^M\in\widehat{\mathcal{U}}$%
			\FOR{$k=0, 1, \ldots$}
			\STATE at all nodes, for $m=1,\ldots,M$, \textbf{do}
		\STATE \ \ \ \ $u_m^{k+1} \coloneqq \mathrm{prox}_{M\omega_m H_m^*/(\gamma_{k} \eta)} \big(u_m^k$ 
		\STATE \ \ \ \ $\ \ \ \ {}+\frac{M\omega_m}{\gamma_k \eta}K_mx_R^k\big)$
\STATE \ \ \ \ $a_m^k \coloneqq M\omega_m x_R^k - K_m^*\big((\gamma_k +\gamma_{k+1})u_m^{k+1}$%
\STATE \ \ \ \ $\ \ \ \ \ {}-\gamma_k u_m^k\big)$
\STATE \ \ \ \ transmit $a_m^{k}$  to master
\STATE at master, \textbf{do}
\STATE \ \ \ \ $x^{k+1}_R \coloneqq \mathrm{prox}_{\gamma_{k+1} R}\big(\frac{1}{M}\sum_{m=1}^M a_m^k\big)$
\STATE \ \ \ \ broadcast $x^{k+1}_R$ to all nodes
			\ENDFOR
		\end{algorithmic}
	\end{algorithm}
\begin{algorithm}[H]
		\caption*{\textbf{Distributed Forward--Backward Alg.}}
		\begin{algorithmic}
			\STATE \textbf{input:} $(\gamma_k)_{k\in\mathbb{N}}$, $x_1\in\mathcal{X}$
			\FOR{$k=1,2, \ldots$}
			 \STATE at all nodes, for $m=1,\ldots,M$, \textbf{do}
			\STATE \ \ \ \ $a_{m}^{k}\coloneqq \nabla F_m(x^{k})$
			\STATE \ \ \ \ transmit $a_m^{k}$  to master
			\STATE at master, \textbf{do}
			 \STATE \ \ \ \ $x^{k+1} \coloneqq \mathrm{prox}_{\gamma_k R}(x^k-\frac{\gamma_{k}}{M}\sum_{m=1}^M  a_m^k)$
			 \STATE \ \ \ \ broadcast $x^{k+1}$ to all nodes
			\ENDFOR
		\end{algorithmic}
	\end{algorithm}
	
	\ 
	\end{minipage}
	\end{figure*}

We now focus on the more general problem \eqref{eqpbdi} and we derive distributed versions of the PD3O and PDDY algorithms to solve it. 
For this, we develop a lifting technique: we recast the minimization of $R(x)+\frac{1}{M}\sum_{m=1}^M \big(F_m(x)+H_m(K_m x)  \big)$ as the minimization of 
\begin{equation*}
\widehat{R}(\hat{x})+\widehat{F}(\hat{x})+\widehat{H}(\widehat{K}\hat{x}),
\end{equation*}
 as follows. Let $(\omega_m)_{m=1}^M$ be a sequence of positive weights, whose sum is $1$; they can be used to mitigate different $\|K_m\|$, by setting $\omega_m \propto 1/\|K_m\|^2$, or different $L_{F_m}$, by setting $\omega_m \propto L_{F_m}^2$,
as a rule of thumb.

 We introduce the Hilbert space $\widehat{\mathcal{X}}=\mathcal{X}\times\cdots\times\mathcal{X}$ ($M$ times), endowed with the inner product 
 \begin{equation*}
 \langle\cdot\,,\cdot\rangle_{\widehat{\mathcal{X}}}:(\hat{x},\hat{x}')\mapsto \sum_{m=1}^M \omega_m \langle x_m,x'_m\rangle,
 \end{equation*}
  and the Hilbert space $\widehat{\mathcal{U}}=\mathcal{U}_1\times\cdots\times\mathcal{U}_M$, endowed with the inner product 
  \begin{equation*}
  \langle\cdot\,,\cdot\rangle_{\widehat{\mathcal{U}}}:(\hat{u},\hat{u}')\mapsto \sum_{m=1}^M \omega_m \langle u_m,u'_m\rangle.
 \end{equation*}
Furthermore, we introduce $\widehat{K}:\hat{x}=(x_m)_{m=1}^M\in\widehat{\mathcal{X}}\mapsto (K_1 x_1,\ldots,K_M x_M)\in\widehat{\mathcal{U}}$, and 
the functions $\imath_=: \hat{x}\in \widehat{\mathcal{X}}\mapsto \{0$ if $x_1=\cdots= x_M,$ $+\infty$ otherwise$\}$, $\widehat{R}:\hat{x}\in\widehat{\mathcal{X}}\mapsto R(x_1)+\imath_=(\hat{x})$, 
 $\widehat{H}:\hat{u}\in\widehat{\mathcal{U}}\mapsto\frac{1}{M}\sum_{m=1}^M H_m(u_m)$, and 
 $\widehat{F}:\hat{x}\in\widehat{\mathcal{X}}\mapsto\frac{1}{M}\sum_{m=1}^M  F_m(x_m)$. We have to be careful when defining the gradient and proximity operators, because of the weighted metrics; see in Section~\ref{secA6} for details. 
 
 Doing these substitutions in the PD3O and PDDY algorithms, we obtain the new Distributed PD3O and Distributed PDDY algorithms, shown above.  Their particular cases, also shown above, are the distributed Davis--Yin algorithm when $K_m\equiv I$ and $\eta=1$, the distributed Loris--Verhoeven algorithm when $R=0$, the distributed Chambolle--Pock algorithm when $F_m\equiv 0$, the distributed Douglas--Rachford algorithm when $F_m\equiv 0$, $K_m\equiv I$ and $\eta=1$, the (classical) distributed forward--backward algorithm when  $H_m\equiv 0$.

We can easily translate Theorems \ref{th1}--\ref{th4} to these distributed algorithms; the corresponding theorems are given in Section~\ref{secA6}. In a nutshell, we obtain the same convergence results and rates with any number of nodes $M\geq 1$ as in the non-distributed setting, 
for any $\gamma_0\in (0,2/L_{\widehat{F}})$ and $\eta\geq \|\widehat{K}\|^2$, where $L_{\widehat{F}}$ and $\widehat{K}$ are detailed in Section~\ref{secA6}.
 Hence, to our knowledge, we are the first to propose distributed proximal splitting methods with guaranteed, possibly accelerated, convergence, to minimize an arbitrary sum of smooth or nonsmooth functions, possibly composed with linear operators.

\section{Experiments}\label{secexp}

\subsection{Image Deblurring Regularized with Total Variation}\label{secexptv}

We first consider the non-distributed problem \eqref{eqpb}, for the imaging inverse problem of deblurring, 
which consists in restoring an image $y$ corrupted by blur and noise \citep{cha16}. So, we set 
\begin{equation*}
F:x\mapsto \frac{1}{2}\|Ax-y\|^2,
\end{equation*}
 where the linear operator $A$ corresponds to a 2-D convolution with a lowpass filter, with $L_F=1$. The filter is approximately Gaussian and chosen so that $F$ is $\mu_F$-strongly convex with $\mu_F=0.01$. $y$ is obtained by applying $A$ to the classical $256\times256$ Shepp--Logan phantom image, with additive Gaussian noise. $R=\imath_{0}$ enforces nonnegativity of the pixel values. $H\circ K$ corresponds to the classical `isotropic' total variation (TV) \citep{cha16,con17}, with $H=0.6$ times the $l_{1,2}$ norm and $K$ the concatenation of vertical and horizontal finite differences. 

We compare the nonaccelerated, i.e.\ with constant $\gamma_k$, and accelerated versions, with decaying $\gamma_k$, of the PD3O, PDDY and Condat--V\~u algorithms. We initialize the dual variables at zero and the estimate of the solution as $y$. We set $\gamma_0=1.7$, $\kappa=0.15$, $\eta=8\geq \|K\|^2$ (except for the accelerated Condat--V\~u algorithm proposed in \cite{cha162}, for which $\eta=16$ and $\gamma=0.5$). 

\begin{figure}[t]
\begin{center}
\begin{tabular}{cc}
\!\!\!\!\!\!\!\!\!\!\!\!\includegraphics[scale=0.62]{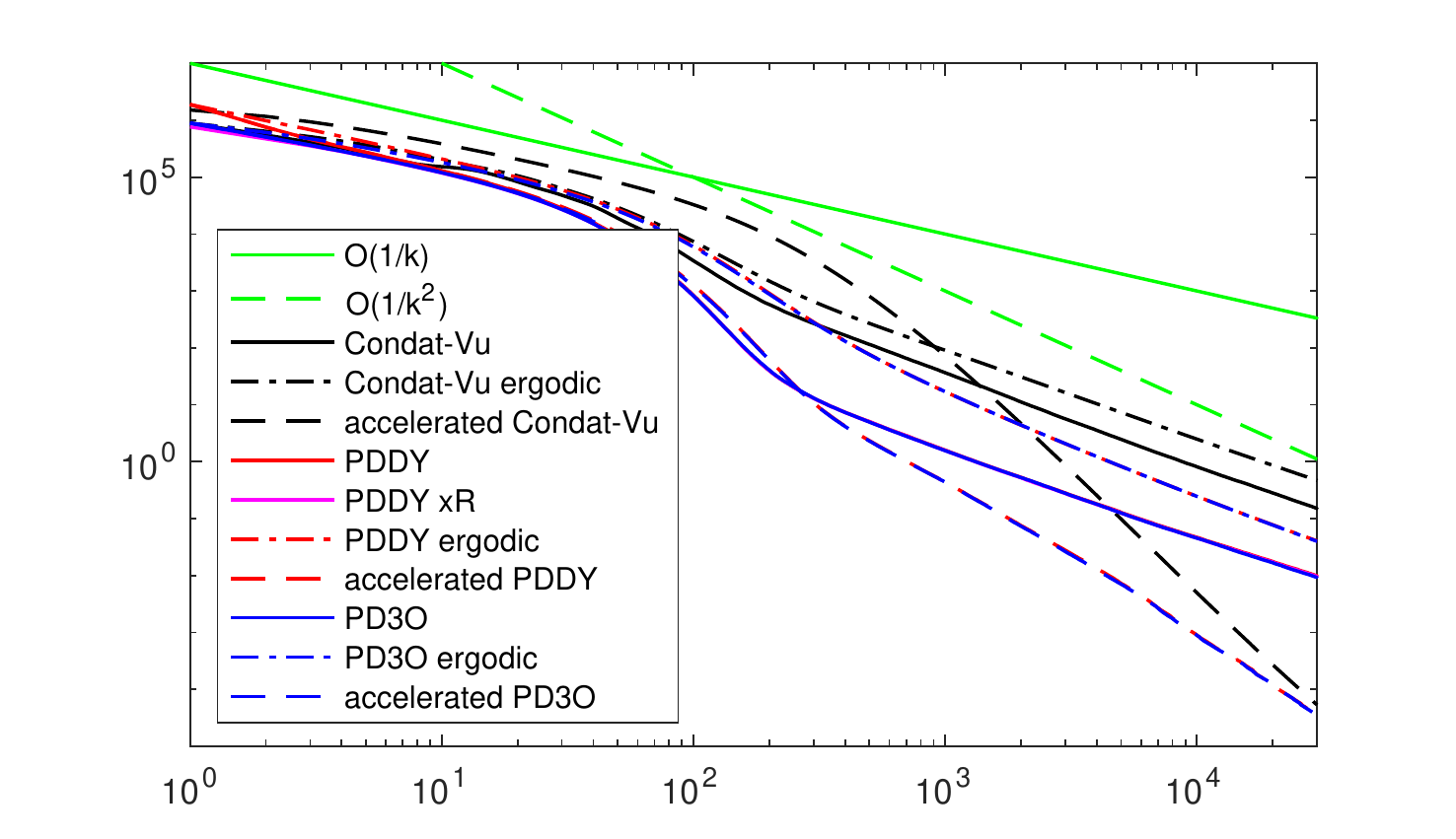}&\!\!\!\!\!\!\!\!\!\!\!\!\!\!\!\!\!\!\includegraphics[scale=0.62]{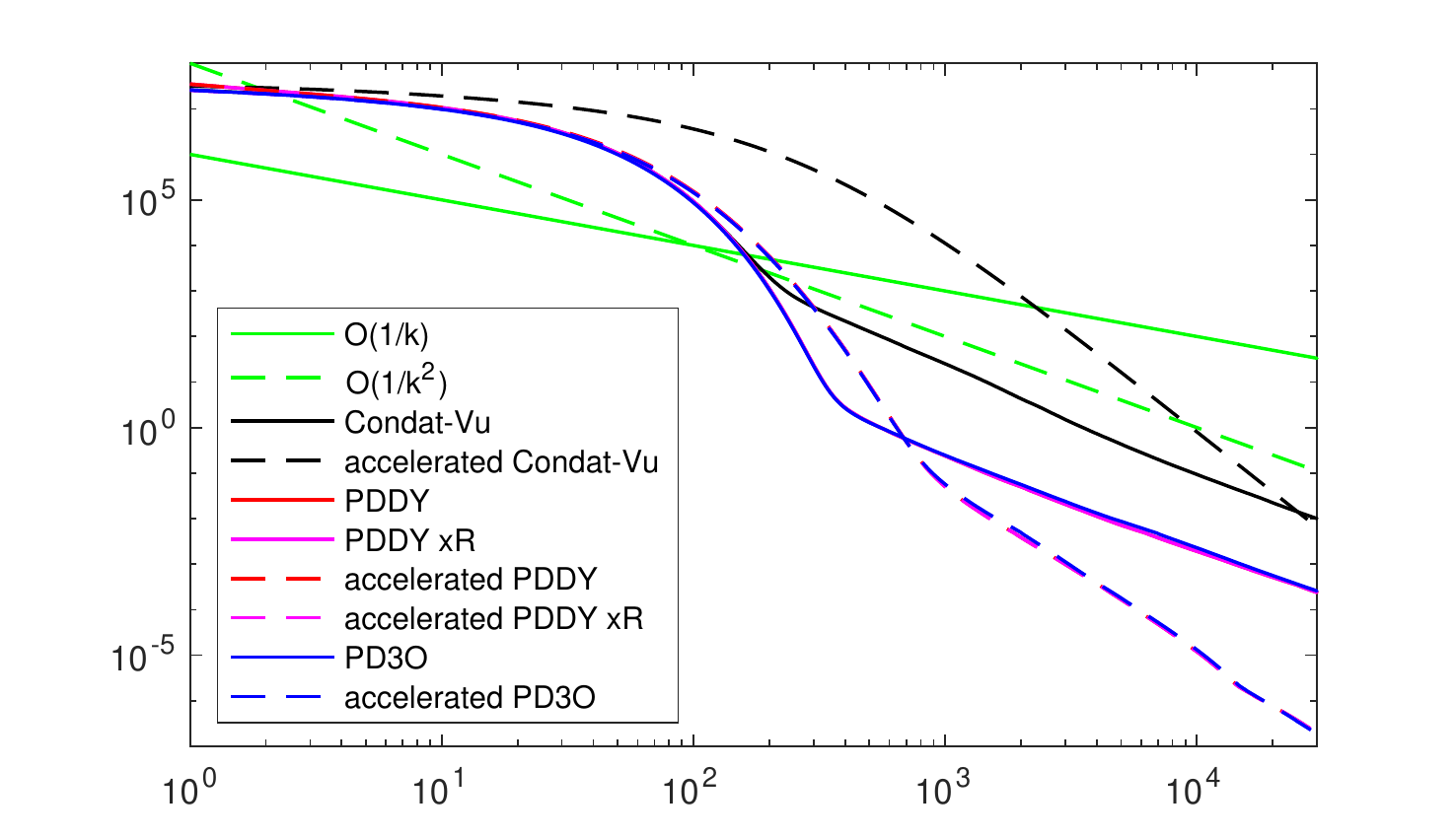}\\
\!\!\!\!\!\!(a) $\Psi(x^k)-\Psi(x^\star)$ w.r.t.\ \# iterations&\!\!\!\!\!\!\!\!(b) $\|x^k-x^\star\|^2$ w.r.t.\ \# iterations
\end{tabular}
\end{center}
\caption{Convergence error, in log-log scale, for the experiment of image deblurring regularized with the total variation, see Section \ref{secexptv} for details.}\label{fig1}
\end{figure}

The results are illustrated in Figure~\ref{fig1} (implementation in Matlab). We observe that the PD3O and PDDY algorithms have almost identical variables: the pink, red, blue curves are superimposed; we know that both algorithms are identical and revert to the Loris--Verhoeven algorithm when $R=0$. Here $R\neq 0$ but the nonnegativity constraint does not change the solution significantly, which explains the similarity of the two algorithms. 

Note that $x^k$ in the PDDY algorithm is not feasible with respect to nonnegativity, and the red curve actually shows $F(x^k)+H(Kx^k)-\Psi(x^\star)$. 
In the nonaccelerated case, $\Psi(x^k)$ decays faster than $O(1/k)$ but slower than $O(1/k^2)$, which is coherent with Theorem~\ref{th1}. The same holds for $\|x^k-x^\star\|^2 \leq \frac{2}{\mu_F}(\Psi(x^k)-\Psi(x^\star))$. 

The accelerated versions improve the convergence speed significantly: $\Psi(x^k)$ and $\|x^k-x^\star\|^2$  decay even faster than $O(1/k^2)$, in line with Theorems \ref{th2} and \ref{th3}. In all cases, the Condat--V\~u algorithm is outperformed. Also, there is no interest in considering the ergodic iterate instead of the last iterate, since the former converges at the same asymptotic rate as the latter, but slower.

\subsection{Image Deblurring Regularized with Huber-TV}\label{secexptv2}

We consider the same deblurring experiment as before, but we make $H$ smooth by taking the Huber function instead of the $l_1$ norm in the total variation; that is, $\lambda |\cdot| $ in the latter is replaced  by 
\begin{equation*}
h:t\in\mathbb{R}\mapsto  \begin{cases} \frac{\lambda}{2\nu}t^2 & \text{if} \quad |t|\leq \nu, \\ 
\lambda \left(|t|-\frac{\nu}{2}\right) & \text{otherwise}, \end{cases}
\end{equation*}
for some $\nu>0$ and $\lambda>0$ (set here as $0.1$ and $0.6$, respectively). We can also write $h$ without branching as  $h(t)=\frac{\lambda}{2\nu}\max(\nu-|t|,0)^2+\lambda(|t|-\frac{\nu}{2})$.
It is known that $h$ is $L_h$-smooth with $L_h=\lambda/\nu$. For any $\gamma>0$ and $t\in\mathbb{R}$, we have 
$\mathrm{prox}_{h^*/\gamma}(t)=t/\max(|t|/\lambda,1+\frac{\nu}{\lambda\gamma} )$. Except for  $H$, everything is unchanged. 

\begin{figure}[t]
\begin{center}
\begin{tabular}{cc}
\!\!\!\!\!\!\!\!\!\!\!\includegraphics[scale=0.615]{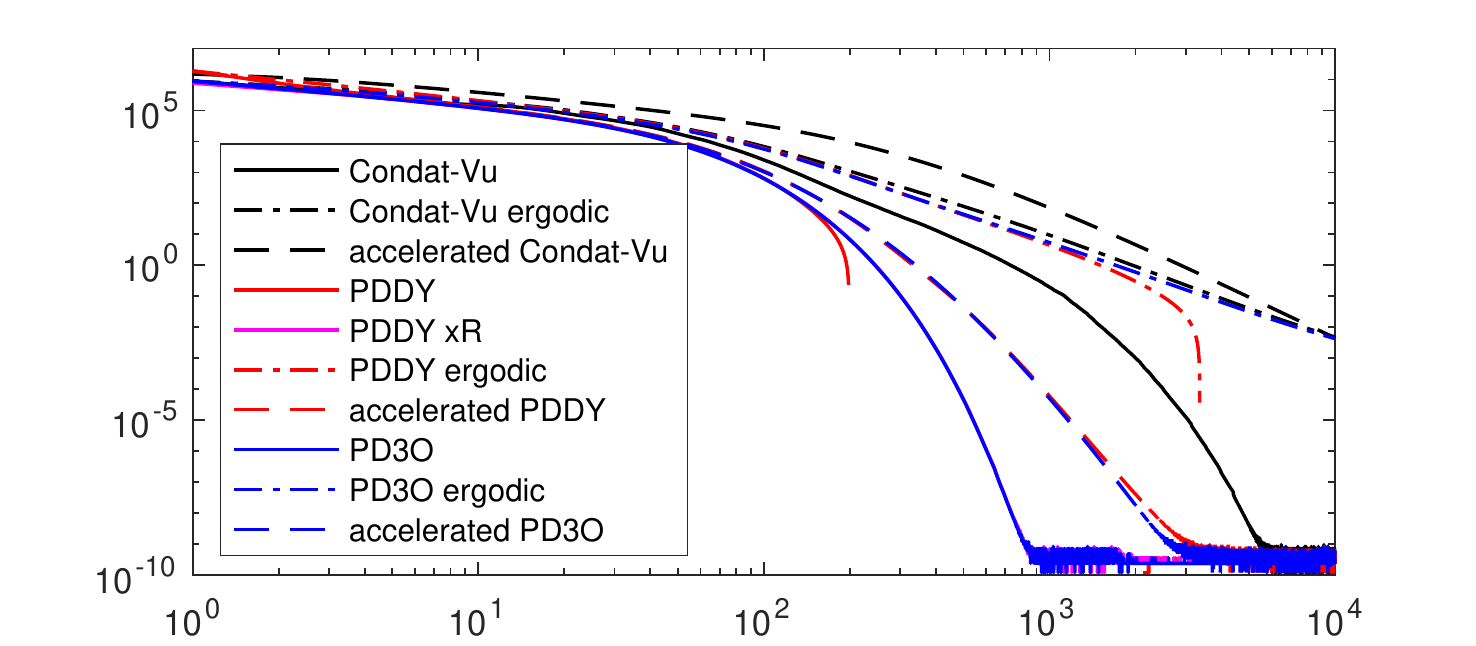}&\!\!\!\!\!\!\!\!\!\!\!\!\!\!\!\!\!\!\includegraphics[scale=0.615]{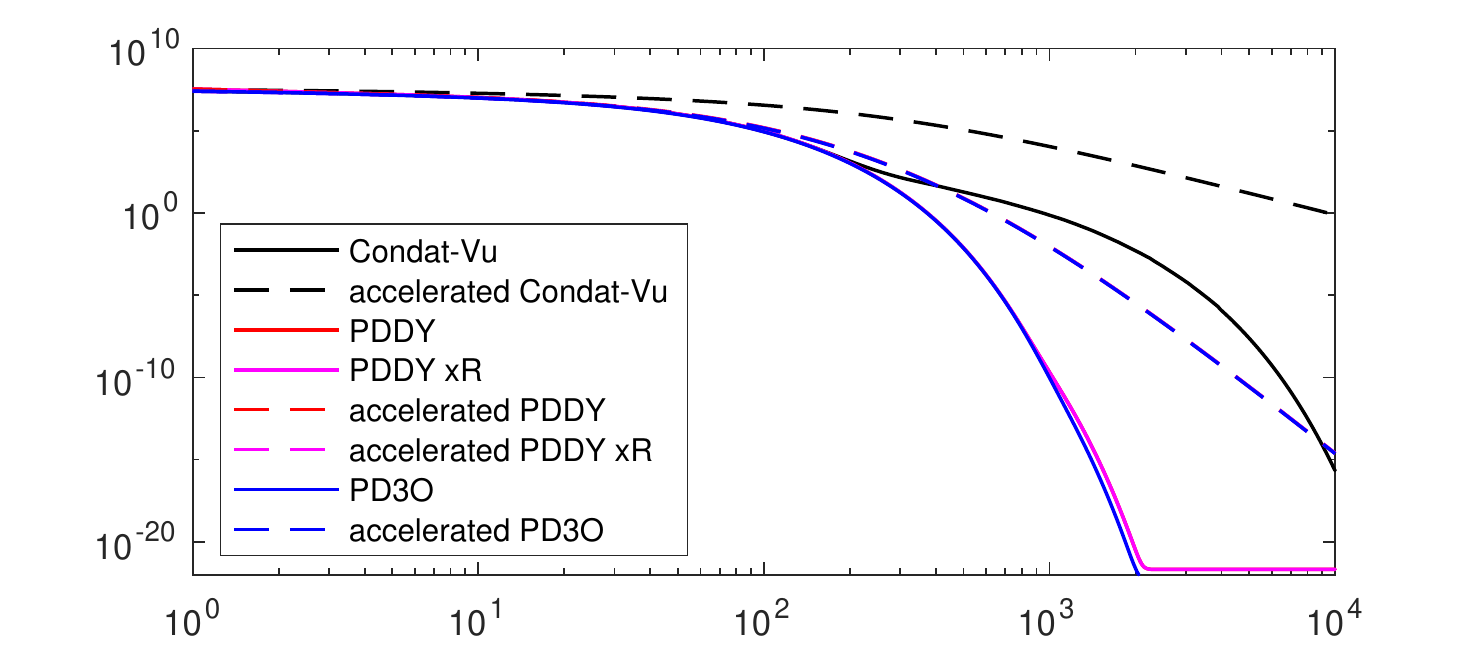}\\
\!\!\!\!\!\!(a) $\Psi(x^k)-\Psi(x^\star)$ w.r.t.\ \#  iterations&\!\!\!\!\!\!\!\!\!\!(b) $\|x^k-x^\star\|^2$ w.r.t.\ \# iterations
\end{tabular}
\end{center}
\caption{Convergence error, in log-log scale, for the experiment of image deblurring regularized with the smooth Huber-total-variation, so that linear convergence occurs, see Section \ref{secexptv2} for details.}\label{fig2}
\end{figure}

The results are illustrated in Figure~\ref{fig2}. Again, 
the PD3O and PDDY algorithms behave very similarly; they converge linearly, as proved in Theorem \ref{th4}, and achieve machine precision in finite time. $x^k$ in the PDDY algorithm is not feasible and $F(x^k)+H(Kx^k)-\Psi(x^\star)$ (red curve) takes negative values (not shown in log scale); so, $x_R^k$ is the variable to study in this setting. We tested the `accelerated' versions of the algorithms with decaying $\gamma_k$, but in this scenario, they are much slower and not suitable.  Again, the Condat--V\~u algorithm is outperformed and the ergodic sequences converge much slower. Interestingly, the image $x^\star$ is visually the same with TV and with Huber-TV.

\subsection{SVM with Hinge Loss}\label{sechl}

Here we consider  Problem \eqref{eqpbdi}  in the special case with $\mathcal{X}=\mathbb{R}^d$, for some $d\geq 1$,
$F_m\equiv 0$, and $K_m\equiv I$; that is, the problem of minimizing
\begin{equation}
\Psi(x) = \frac{1}{M}\sum_{m=1}^M H_m(x) +  R(x).\label{eq:no98f08hf}
\end{equation}
In particular, to train a binary classifier, we consider the classical SVM problem with hinge loss, 
which has the form \eqref{eq:no98f08hf} with
 $R(x)=\frac{\alpha}{2}\|x\|^2$, for some $\alpha>0$, and $H_m(x)=\max(1-b_m a_m^\mathrm{T}x,0)$, with data samples $a_m\in\mathbb{R}^d$ and $b_m\in\{-1,1\}$.

For any $\gamma>0$ we have $\mathrm{prox}_{\gamma R}(x)= x/(1+\gamma\alpha)$.
We could view the dot product $x\mapsto  b_m a_m^\mathrm{T}x$ as a linear operator $K_m$, but it is more interesting to integrate it in the function $H_m$. Indeed, as is perhaps not well known, the proximity operator of $H_m$ has a closed form: for any $\gamma>0$,
\begin{equation*}
\mathrm{prox}_{\gamma H_m}:x\in\mathbb{R}^d \mapsto x-\frac{b_m}{\eta_m}\max\big(\min(b_m a_m^\mathrm{T}x-1,0),-\eta_m\gamma\big)a_m,
\end{equation*}
 where $\eta_m= a_m^\mathrm{T}a_m=\|a_m\|^2$.
Thus, we use the Distributed Douglas--Rachford algorithm, a particular case of the distributed PD3O and PDDY algorithms. Since $R$ is $\alpha$-strongly convex, we also use the accelerated version of the algorithm with varying stepsizes, like in Theorem~\ref{th2}. We can note that 
 in the context of Federated learning \citep{kon16,mal20}, where each $m$ corresponds to the smart phone or computer of a  different user with its own data $(a_m,b_m)$ stored locally, the problem is solved in a collaborative way  but with preserved privacy, without the users sharing their data.

\begin{figure}[t]
\begin{center}
\begin{tabular}{cc}
\hspace{-5mm}
\!\!\includegraphics[scale=0.60]{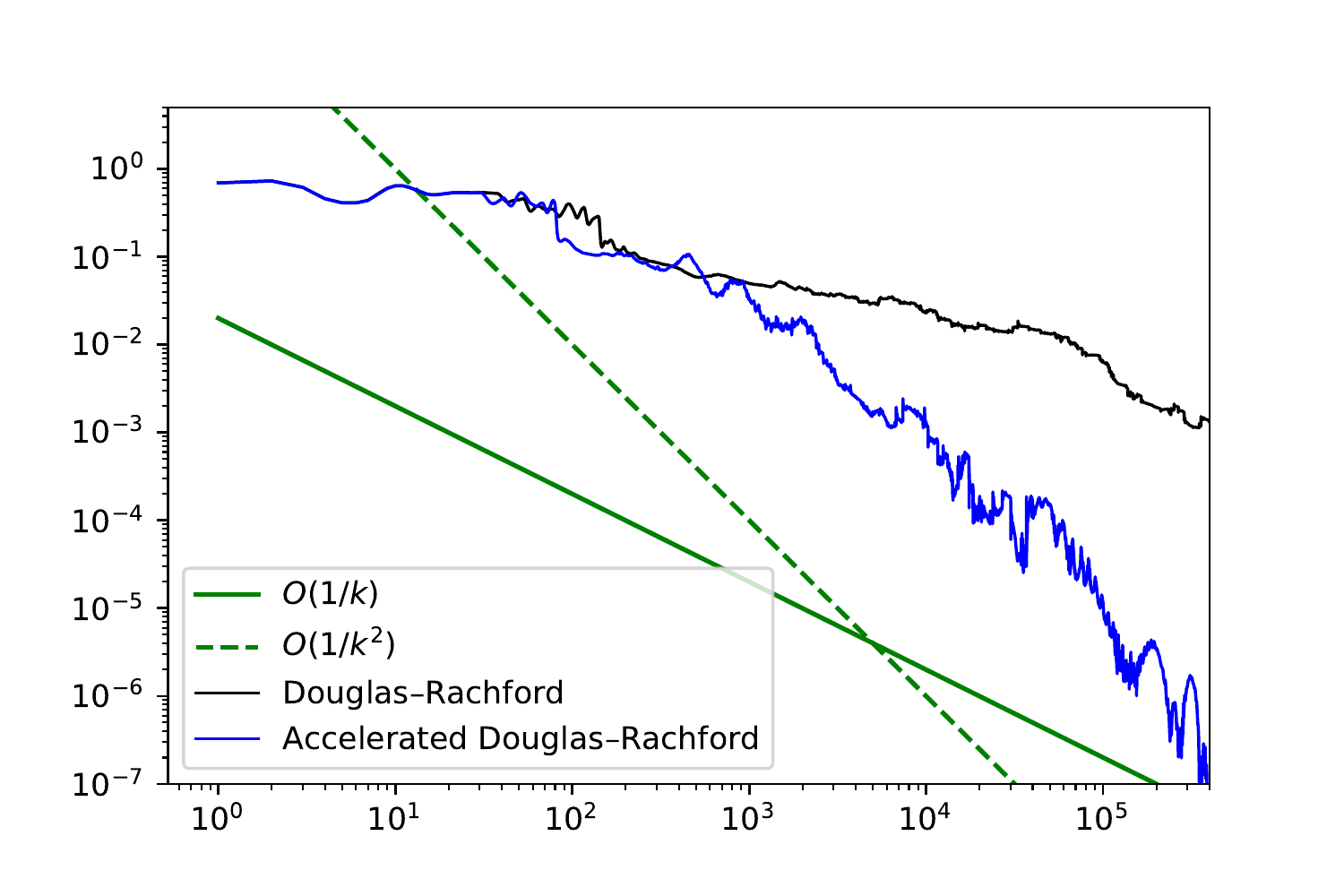}&\hspace{-9mm}\!\!\!\!
\includegraphics[scale=0.60]{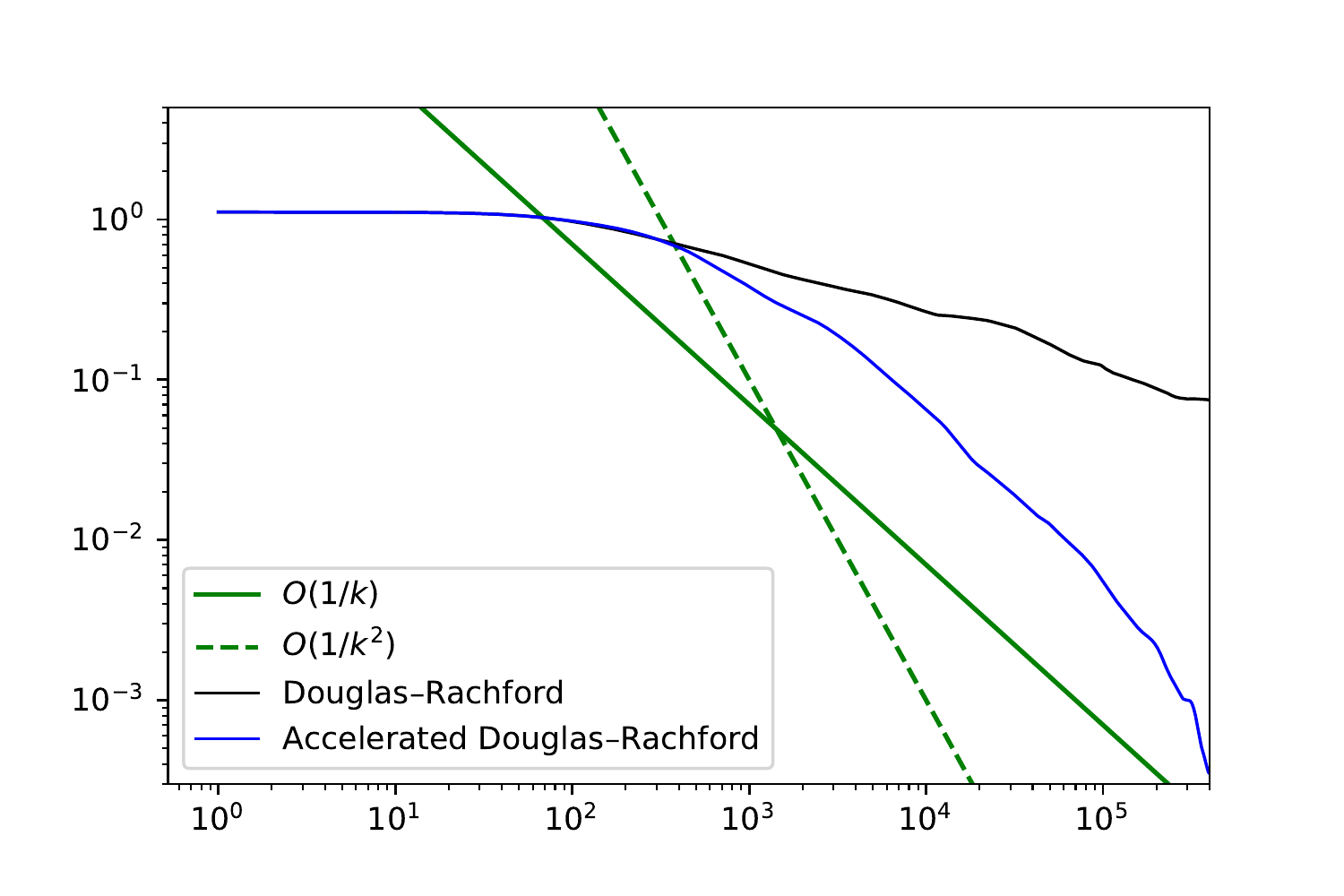}\hspace{-5mm}\\
\!\!\!\!\!\!(a) $\Psi(x^k)-\Psi(x^\star)$ w.r.t.\ \# iterations&\!\!\!\!\!\!(b) $\|x^k-x^\star\|^2$ w.r.t.\ \# iterations
\end{tabular}
\end{center}
\caption{Convergence error, in log-log scale, for the SVM binary classification experiment with hinge loss, see Section \ref{sechl} for details.}\label{fig3}
\end{figure}

The method was implemented in Python on a single machine and tested on the dataset `australian'  from  the LibSVM base~\citep{chang2011libsvm}, with $d=15$ and $M=680$. 
We set $\omega_m \equiv 1/M$, $\alpha=0.1$, $\gamma_0 = 0.1$, and we used zero vectors for the initialization. The results are shown in Figure~\ref{fig3}. Despite the oscillations, we observe that both the objective suboptimality and the squared distance to the solution converge sublinearly, with rates looking like $o(1/\sqrt{k})$ and $O(1/k^2)$ for the nonaccelerated and accelerated algorithms, respectively,
as guaranteed by Theorems \ref{th1} and \ref{th2}. The proposed accelerated version of the distributed Douglas--Rachford algorithm yields a significant speedup.

\section{Derivation of the Algorithms}\label{secA5}

In this section, we give the details of the derivation of the PD3O and PPDY algorithms, and their particular cases, to solve:
\begin{equation*}
{\minimize_{x\in\mathcal{X}} \; F(x)+R(x)+H(K x)},
 \end{equation*}
with same notations and assumptions as above. Let $\eta\geq \|K\|^2$,  let $\mathcal{W}$ be a real Hilbert space and $C:\mathcal{W}\rightarrow \mathcal{U}$ be a linear operator, such that $KK^*+CC^*=\eta I$.  We set $Q:(x,w)\mapsto Kx+Cw$. We have $QQ^*=\eta I$. Let $(\gamma_k)_{k\in\mathbb{N}}$ be a sequence of positive stepsizes. 

\subsection{The Davis--Yin Algorithm}

In this section, we state the results on the Davis--Yin algorithm, which we be needed to analyze the PD3O and PPDY algorithms.

The Davis--Yin algorithm to minimize the sum of 3 convex functions $\widetilde{F}+G+J$ over a real Hilbert space $\mathcal{Z}$ (assuming that there exists a solution $z^\star$ such that $0\in \nabla\widetilde{F}(z^\star)+\partial G(z^\star)+\partial J(z^\star)$) is~\citep{dav17}:

Let $z_J^0\in\mathcal{Z}$, $u_G^0\in\mathcal{Z}$. For $k=0,1,\ldots$ iterate:
\begin{equation}\left\lfloor\begin{array}{l}
z^{k+1}_G = \mathrm{prox}_{\gamma_k G}(z_J^k + \gamma_k u_G^k)\\
u_G^{k+1} = u_G^k+\frac{1}{\gamma_k}(z_J^k  - z_G^{k+1} )\\
z^{k+1}_J = \mathrm{prox}_{\gamma_{k+1} J}\big(z_G^{k+1} - \gamma_{k+1} u_G^{k+1}- \gamma_{k+1}\nabla \widetilde{F} (z^{k+1}_G)\big).
\end{array}\right.\label{eqdybasea1}\end{equation}
Equivalently, introducing the variable $r^k\coloneqq z_J^k + \gamma_k u_G^k$:
 let $r^0\in\mathcal{Z}$.  
For $k=0,1,\ldots$ iterate:
\begin{equation}\left\lfloor\begin{array}{l}
z^{k+1}_G = \mathrm{prox}_{\gamma_k G}(r^k)\\
z^{k+1}_J = \mathrm{prox}_{\gamma_{k+1} J}\big((1+ \frac{\gamma_{k+1}}{\gamma_k})z_G^{k+1} - \frac{\gamma_{k+1}}{\gamma_k}r^k
- \gamma_{k+1}\nabla \widetilde{F} (z^{k+1}_G)\big)\\
r^{k+1}=z^{k+1}_J + \frac{\gamma_{k+1}}{\gamma_k}(r^k - z_G^{k+1} ).
\end{array}\right.\label{eqdy2}\end{equation}
Equivalently: let $r^0\in\mathcal{Z}$.  
For $k=0,1,\ldots$ iterate:
\begin{equation}\left\lfloor\begin{array}{l}
z^{k+1}_G = \mathrm{prox}_{\gamma_k G}(r^k)\\
u^{k+1}_J = \mathrm{prox}_{J^*/\gamma_{k+1} }\big((\frac{1}{\gamma_{k+1}}+ \frac{1}{\gamma_k})z_G^{k+1} - \frac{1}{\gamma_k}r^k - \nabla \widetilde{F} (z^{k+1}_G)\big)\\
r^{k+1}=z_G^{k+1} 
- \gamma_{k+1}\nabla \widetilde{F} (z^{k+1}_G)-\gamma_{k+1}u^{k+1}_J.
\end{array}\right.\label{dyapp1}\end{equation}

In our notations, Theorem 3.3 of \cite{dav17} translates into Lemma~\ref{lem1} as follows; we assume that $\widetilde{F}$ is $L_{\widetilde{F}}$-smooth and $\mu_{\widetilde{F}}$-strongly convex and that $G$ is $\mu_G$-strongly convex, for some  $L_{\widetilde{F}}>0$, $\mu_{\widetilde{F}}\geq 0$, $\mu_G\geq 0$.

\begin{lemma}[accelerated Davis--Yin algorithm]\label{lem1}Suppose that $\mu_{\widetilde{F}}+\mu_G>0$. Let $z^\star$ be the unique minimizer of $\widetilde{F}+G+J$; that is, $0\in \nabla\widetilde{F}(z^\star)+\partial G(z^\star)+\partial J(z^\star)$. Let $u_G^\star$ be such that $u_G^\star\in\partial G(z^\star)$ and $0\in \nabla\widetilde{F}(z^\star)+ \partial J(z^\star)+u_G^\star$. Let $\kappa\in(0,1)$ and $\gamma_0\in (0,2(1-\kappa)/L_{\widetilde{F}})$.  Set $\gamma_1=\gamma_0$ and
\begin{equation*}
\gamma_{k+1}=\frac{-\gamma_k^2\mu_{\widetilde{F}}\kappa+\gamma_k\sqrt{(\gamma_k\mu_{\widetilde{F}}\kappa)^2+1+2\gamma_k\mu_G}}{1+2\gamma_k\mu_G}
,\quad\mbox{for every }k\geq 1.
\end{equation*}
Then, for every $k\geq 1$,
\begin{equation*}
\|z_G^{k+1}-z^\star\|^2\leq \frac{\gamma_{k+1}^2}{1-\gamma_{k+1}\mu_{\widetilde{F}}\kappa}c_0
=O\big(1/k^2\big),
\end{equation*}
where
\begin{equation*}
c_0=\frac{1-\gamma_0\mu_{\widetilde{F}}\kappa}{\gamma_0^2}\|z_G^1-z^\star\|^2+\|
u_G^1-u_G^\star\|^2.
\end{equation*}
Note that $u_G^1=(r^0-z_G^1)/\gamma_0$.
\end{lemma}

Linear convergence occurs in the following conditions, according to Theorem D.6 in the preprint version of \cite{dav17}, which translates into Lemma~\ref{lem2} as follows.
We assume that $\widetilde{F}$ is $L_{\widetilde{F}}$-smooth and $\mu_{\widetilde{F}}$-strongly convex, $G$ is $\mu_G$-strongly convex, and $J$ is $\mu_J$-strongly convex, for some  $L_{\widetilde{F}}>0$, $\mu_{\widetilde{F}}\geq 0$, $\mu_G\geq 0$, $\mu_J\geq 0$. We consider constant stepsizes $\gamma_k\equiv \gamma$, for some $\gamma\in (0,2/L_{\widetilde{F}})$.

\begin{lemma}[linear convergence of the Davis--Yin algorithm]\label{lem2}Suppose that $\mu_{\widetilde{F}}+\mu_G+\mu_J>0$ and that $G$ is $L_G$-smooth, for some $L_G>0$, or $J$ is $L_J$-smooth, for some $L_J>0$. Let $z^\star$ be the unique minimizer of $\widetilde{F}+G+J$; that is, $0\in \nabla\widetilde{F}(z^\star)+\partial G(z^\star)+\partial J(z^\star)$. The dual problem of minimizing $(\widetilde{F}+J)^*(-u)+G^*(u)$ over $u\in\mathcal{Z}$ is strongly convex too; let $u_G^\star$ be its unique solution. We have  $u_G^\star\in\partial G(z^\star)$ and $0\in \nabla\widetilde{F}(z^\star)+ \partial J(z^\star)+u_G^\star$. Set $r^\star=z^\star+\gamma u_G^\star$. Then, the Davis--Yin algorithm \eqref{eqdy2} converges linearly: there exists $\rho\in(0,1]$ such that, for every $k\in \mathbb{N}$,
\begin{equation}
\|r^k-r^\star\|^2\leq (1-\rho)^k \|r^0-r^\star\|^2.\label{eqxx0}
\end{equation}
\end{lemma}

Loose lower bounds for $\rho$ are given in \citealp[Theorem D.6]{dav17}. \medskip

We have the following corollary of Lemma~\ref{lem2}:

\begin{proposition}[linear convergence of the other variables in the Davis--Yin algorithm]\label{prop1}In the same conditions and notations as in Lemma~\ref{lem2}, we have, for every $k\in \mathbb{N}$, 
\begin{align}
\|z_G^{k+1}-z^\star\|^2&\leq (1-\rho)^k \|r^0-r^\star\|^2\label{eqxx1}\\
 \|z_J^{k+1}-z^\star\|^2&\leq 4 (1-\rho)^k \|r^0-r^\star\|^2\notag.
 \end{align}
 Also, in the form \eqref{dyapp1} of the algorithm, 
   \begin{equation*}
 \|u_J^{k+1}+u_G^\star+\nabla\widetilde{F}(z^\star)\|^2 \leq \frac{4}{\gamma^2} (1-\rho)^k \|r^0-r^\star\|^2
 \end{equation*}
and, in the form \eqref{eqdybasea1} of the algorithm, 
   \begin{equation*}
\|u_G^{k+1} -u_G^\star\|^2 \leq \frac{1}{\gamma^2}(1-\rho)^k\|r^0-r^\star\|^2.
\end{equation*}
 \end{proposition}
\proof
Let $k\in\mathbb{N}$. 
 By nonexpansiveness of the proximity operator, in view of the first line in \eqref{eqdy2}, we have $\|z^{k+1}_G -z^\star\| \leq \|r^k-r^\star\|$, so that \eqref{eqxx1} follows from \eqref{eqxx0}. In addition, in view of the second line in \eqref{eqdy2}, we have
  \begin{align*}
  \|z_J^{k+1}-z^\star\|^2&\leq  \|2(z_G^{k+1}-z^\star)-(r^k-r^\star)-\gamma(\nabla\widetilde{F}(z_G^{k+1})-\nabla\widetilde{F}(z^\star))\|^2\\
&= \|(z_G^{k+1}-z^\star)-(r^k-r^\star)+(I-\gamma \nabla\widetilde{F})(z_G^{k+1})-(I-\gamma \nabla\widetilde{F})(z^\star)\|^2\\
&= \|(I-\mathrm{prox}_{\gamma G})(r^k)-(I-\mathrm{prox}_{\gamma G})(r^\star)+(I-\gamma \nabla\widetilde{F})(z_G^{k+1})-(I-\gamma \nabla\widetilde{F})(z^\star)\|^2
\end{align*}
and, by nonexpansiveness of $I-\mathrm{prox}_{\gamma G}$ and $I-\gamma \nabla\widetilde{F}$, 
\begin{align*}
  \|z_J^{k+1}-z^\star\|^2&\leq \big(\|r^k-r^\star\|+\|z_G^{k+1}-z^\star\|\big)^2\\
  &\leq 4  \|r^k-r^\star\|^2.
\end{align*}
Using the same arguments, in view of the second line in \eqref{dyapp1}, 
\begin{align*}
  \|u_J^{k+1}+u_G^\star+\nabla\widetilde{F}(z^\star)\|^2 &\leq
     \frac{1}{\gamma^2}\big(\|r^k-r^\star\|+\|z_G^{k+1}-z^\star\|\big)^2\\
    &\leq \frac{4}{\gamma^2}  \|r^k-r^\star\|^2.
\end{align*}
Finally, as visible in the first line of \eqref{eqdy100}, 
since  $r^k= z_J^k + \gamma_k u_G^k$, and using the Moreau identity, we have $u_G^{k+1} = \mathrm{prox}_{G^*/\gamma}(\frac{1}{\gamma}z_J^k +  u_G^k)=
 \mathrm{prox}_{G^*/\gamma}(\frac{1}{\gamma}r^k)$, so that 
\begin{equation*}
\|u_G^{k+1} -u_G^\star\|^2 \leq \frac{1}{\gamma^2}\|r^k-r^\star\|^2 .
\end{equation*}
  \endproof

 \subsection{The PD3O Algorithm}

We set $\mathcal{Z}=\mathcal{X}\times\mathcal{W}$, $\widetilde{F}$, 
$G=\widetilde{R}$, $J=\widetilde{H}$, as defined in Section~\ref{secm1}. Doing the substitutions in \eqref{dyapp1}, we get the algorithm:

Let $s^0\in\mathcal{X}$ and $r_w^{0}\in\mathcal{W}$. 
For $k=0,1,\ldots$ iterate:
\begin{equation*}\left\lfloor\begin{array}{l}
x^{k+1} = \mathrm{prox}_{\gamma_k R}(s^k)\\
u^{k+1} = \mathrm{prox}_{H^*/(\gamma_{k+1} \eta)} \Big(K\big((\frac{1}{\gamma_{k+1}}+ \frac{1}{\gamma_k})x^{k+1} - \frac{1}{\gamma_k}s^k - \nabla F(x^{k+1})\big)/\eta-Cr_w^k/(\gamma_k\eta)\Big)\\
s^{k+1}=x^{k+1} 
- \gamma_{k+1}\nabla F (x^{k+1})-\gamma_{k+1}K^*u^{k+1} \\
r_w^{k+1}=-\gamma_{k+1}C^*u^{k+1}.
\end{array}\right.\end{equation*}
We can remove the variable $r_w$ and the algorithm becomes: 
Let $s^{0}\in\mathcal{X}$ and $u^{0}\in\mathcal{U}$. 
For $k=0,1,\ldots$ iterate:
\begin{equation*}\left\lfloor\begin{array}{l}
x^{k+1} = \mathrm{prox}_{\gamma_k R}(s^k)\\
u^{k+1} = \mathrm{prox}_{H^*/(\gamma_{k+1} \eta)} \Big(\frac{1}{\eta}K\big((\frac{1}{\gamma_{k+1}}+ \frac{1}{\gamma_k})x^{k+1} - \frac{1}{\gamma_k}s^k - \nabla F(x^{k+1})\big)+\frac{1}{\eta}CC^*u^k\Big)\\
s^{k+1}=x^{k+1} 
- \gamma_{k+1}\nabla F (x^{k+1})-\gamma_{k+1}K^*u^{k+1} .
\end{array}\right.\end{equation*}
After replacing $CC^*$ by $\eta I-KK^*$, the iteration becomes:
\begin{equation*}\left\lfloor\begin{array}{l}
x^{k+1} = \mathrm{prox}_{\gamma_k R}(s^k)\\
u^{k+1} = \mathrm{prox}_{H^*/(\gamma_{k+1} \eta)} \Big(u^k+\frac{1}{\eta}K\big((\frac{1}{\gamma_{k+1}}+ \frac{1}{\gamma_k})x^{k+1} - \frac{1}{\gamma_k}s^k - \nabla F(x^{k+1})-K^* u^k\big)\Big)\\
s^{k+1}=x^{k+1} 
- \gamma_{k+1}\nabla F (x^{k+1})-\gamma_{k+1}K^*u^{k+1}.
\end{array}\right.\end{equation*}
We can change the variables, so that only one call to $\nabla F$ and $K^*$ appears, which yields the algorithm: 
Let $q^{0}\in\mathcal{X}$ and $u^{0}\in\mathcal{U}$. For $k=0,1,\ldots$ iterate:
\begin{equation*}\left\lfloor\begin{array}{l}
x^{k+1} = \mathrm{prox}_{\gamma_k R}\big(\gamma_{k}(q^k-K^*u^{k}) \big)\\
q^{k+1}=\frac{1}{\gamma_{k+1}}x^{k+1} 
- \nabla F (x^{k+1})\\
u^{k+1} = \mathrm{prox}_{H^*/(\gamma_{k+1} \eta)} \big(u^k+\frac{1}{\eta}K(\frac{1}{\gamma_k}x^{k+1}+q^{k+1}-q^k )\big).
\end{array}\right.\end{equation*}

When $\gamma_k\equiv \gamma$ is constant, we recover the PD3O algorithm~\citep{yan18}.

To derive Theorem~\ref{th2} from Lemma~\ref{lem1}, we simply have to notice that the variable $z_G^{k+1}$ in the latter corresponds to the pair $(x^{k+1},0)$. Also, in the conditions of Theorem~\ref{th2}, let $u^\star$ be any solution of \eqref{eqpd}; that is, $u^\star\in\partial H(Kx^\star)$ and $0\in \partial R(x^\star)+\nabla F(x^\star)+K^* u^\star$. Then the constant $c_0$ is
\begin{equation*}
c_0=\frac{1-\gamma_0\mu_F\kappa}{\gamma_0^2}\|x^1-x^\star\|^2+\|
q^0-{\textstyle\frac{1}{\gamma_0}}x^1-K^*(u^0-u^\star)
+\nabla F(x^\star)\|^2+\eta \|u^0-u^\star\|^2-\|K^*(u^0-u^\star)\|^2.
\end{equation*}

If $K=I$ and $\eta=1$, the PD3O algorithm reverts to the Davis--Yin algorithm, as given in \eqref{eqdybase}. In the conditions of Theorem~\ref{th2}, let $u^\star$ be any solution of \eqref{eqpd}; that is, $u^\star\in\partial H(x^\star)$ and $0\in \partial R(x^\star)+\nabla F(x^\star)+u^\star$. Then the constant $c_0$ is
\begin{equation}
c_0=\frac{1-\gamma_0\mu_F\kappa}{\gamma_0^2}\|x^1-x^\star\|^2+\|
{\textstyle\frac{1}{\gamma_0}}(s^0-x^1)+u^\star+\nabla F(x^\star)\|^2.\label{eqc01}
\end{equation}

\subsection{The PDDY Algorithm}

The PDDY algorithm is obtained like the PD3O algorithm from the David--Yin algorithm, but after swapping the roles of $\widetilde{H}$ and $\widetilde{R}$.

To obtain the PDDY algorithm, starting from \eqref{eqdybasea1}, let us first write the Davis--Yin algorithm as:
Let $z_J^{0}\in\mathcal{Z}$ and $u_G^{0}\in\mathcal{Z}$. 
For $k=0,1,\ldots$ iterate:
\begin{equation*}\left\lfloor\begin{array}{l}
u_G^{k+1} = \mathrm{prox}_{G^*/\gamma_k}(\frac{1}{\gamma_k}z_J^k +  u_G^k)\\
z_G^{k+1}=z_J^k  - \gamma_k(u_G^{k+1}-u_G^{k})\\
z^{k+1}_J= \mathrm{prox}_{\gamma_{k+1} J}\big(z_G^{k+1} -\gamma_{k+1}\nabla \widetilde{F}(z^{k+1}_G)- \gamma_{k+1} u_G^{k+1}\big).
\end{array}\right.\end{equation*}
Equivalently:
Let $r^{0}\in\mathcal{Z}$.  For $k=0,1,\ldots$ iterate:
\begin{equation}\left\lfloor\begin{array}{l}
u_G^{k+1} = \mathrm{prox}_{ G^*/\gamma_k}(r^k/\gamma_k )\\
z_G^{k+1}=r^k  - \gamma_k u_G^{k+1}\\
z^{k+1}_J = \mathrm{prox}_{\gamma_{k+1} J}\big(z_G^{k+1} -\gamma_{k+1}\nabla \widetilde{F}(z^{k+1}_G)- \gamma_{k+1} u_G^{k+1}\big)\\
r^{k+1} = z^{k+1}_J + \gamma_{k+1} u_G^{k+1}.
\end{array}\right.\label{eqdy100}\end{equation}

We set $\mathcal{Z}=\mathcal{X}\times\mathcal{W}$, $\widetilde{F}$, $G=\widetilde{H}$, $J=\widetilde{R}$, as defined in Section~\ref{secm1}. Doing the substitutions in \eqref{eqdy100}, we get the algorithm: 
Let $r_x^{0}\in\mathcal{X}$, $r_w^{0}\in\mathcal{W}$. For $k=0,1,\ldots$ iterate:
\begin{equation*}\left\lfloor\begin{array}{l}
u^{k+1} = \mathrm{prox}_{H^*/(\gamma_{k} \eta)} \big((Kr_x^k+Cr_w^k)/(\gamma_k \eta)\big)\\
x^{k+1}=r_x^k  - \gamma_k K^*u^{k+1}\\
x^{k+1}_R = \mathrm{prox}_{\gamma_{k+1} R}\big(x^{k+1} -\gamma_{k+1}\nabla F (x^{k+1})- \gamma_{k+1} K^*u^{k+1}\big)\\
r_x^{k+1} = x^{k+1}_R + \gamma_{k+1} K^*u^{k+1}\\
r_w^{k+1} =  \gamma_{k+1} C^*u^{k+1}.
\end{array}\right.\end{equation*}
We can remove the variable $r_w$ and rename $r_x$ as $s$:
\begin{equation*}\left\lfloor\begin{array}{l}
u^{k+1} = \mathrm{prox}_{H^*/(\gamma_{k} \eta)} \big(Ks^k/(\gamma_k \eta)+CC^*u^k/ \eta\big)\\
x^{k+1}=s^k  - \gamma_k K^*u^{k+1}\\
x^{k+1}_R = \mathrm{prox}_{\gamma_{k+1} R}\big(x^{k+1} -\gamma_{k+1}\nabla F (x^{k+1})- \gamma_{k+1} K^*u^{k+1}\big)\\
s^{k+1} = x^{k+1}_R + \gamma_{k+1} K^*u^{k+1}.
\end{array}\right.\end{equation*}
The algorithm becomes: 
Let $s^{0}\in\mathcal{X}$, $u^{0}\in\mathcal{U}$. For $k=0,1,\ldots$ iterate:
\begin{equation*}\left\lfloor\begin{array}{l}
u^{k+1} = \mathrm{prox}_{H^*/(\gamma_{k} \eta)} \big(u^k +K(s^k/\gamma_k - K^*u^k)/ \eta\big)\\
x^{k+1}=s^k  - \gamma_k K^*u^{k+1}\\
x^{k+1}_R = \mathrm{prox}_{\gamma_{k+1} R}\big(x^{k+1} -\gamma_{k+1}\nabla F (x^{k+1})- \gamma_{k+1} K^*u^{k+1}\big)\\
s^{k+1} = x^{k+1}_R + \gamma_{k+1} K^*u^{k+1}.
\end{array}\right.\end{equation*}
Equivalently: 
Let $x_R^{0}\in\mathcal{X}$, $u^{0}\in\mathcal{U}$. For $k=0,1,\ldots$ iterate:
\begin{equation*}\left\lfloor\begin{array}{l}
u^{k+1} = \mathrm{prox}_{H^*/(\gamma_{k} \eta)} \big(u^k +Kx_R^k/(\gamma_k \eta)\big)\\
x^{k+1}=x_R^k  - \gamma_k K^*(u^{k+1}-u^k)\\
x^{k+1}_R= \mathrm{prox}_{\gamma_{k+1} R}\big(x^{k+1} -\gamma_{k+1}\nabla F (x^{k+1})- \gamma_{k+1} K^*u^{k+1}\big).
\end{array}\right.\end{equation*}
We can write the algorithm with only one call of $K^*$ per iteration by introducing an additional variable $p$: 
Let $x_R^{0}\in\mathcal{X}$, $u^{0}\in\mathcal{U}$. Set $p^{0}=K^*u^{0}$. For $k=0,1,\ldots$ iterate:
\begin{equation*}\left\lfloor\begin{array}{l}
u^{k+1} = \mathrm{prox}_{H^*/(\gamma_{k} \eta)} \big(u^k +\frac{1}{\gamma_k \eta}Kx_R^k\big)\\
p^{k+1} = K^*u^{k+1}\\
x^{k+1}=x_R^k  - \gamma_k (p^{k+1}-p^k)\\
x^{k+1}_R = \mathrm{prox}_{\gamma_{k+1} R}\big(x^{k+1} -\gamma_{k+1}\nabla F (x^{k+1})- \gamma_{k+1} p^{k+1}\big).
\end{array}\right.\end{equation*}

When $\gamma_k\equiv \gamma$ is constant, we recover the PDDY algorithm~\citep{sal20}.

Let us now derive Theorem~\ref{th3} from Lemma~\ref{lem1}. The variable $z_G^{k+1}$ in the latter corresponds to  the pair $\big(x^{k+1},\gamma_k C^*(u^k-u^{k+1})\big)$, so that  $\|z_G^{k+1}-z^\star\|^2$ becomes 
\begin{align}
\|x^{k+1}-x^\star\|^2 + \|\gamma_k C^*(u^k-u^{k+1})\|^2&=
\|x^{k+1}-x^\star\|^2 + \gamma_k^2 \langle CC^* (u^k-u^{k+1}),u^k-u^{k+1}\rangle\notag\\
&=
\|x^{k+1}-x^\star\|^2 + \gamma_k^2 \langle (\eta I - KK^*) (u^k-u^{k+1}),u^k-u^{k+1}\rangle\notag\\
&=
\|x^{k+1}-x^\star\|^2 + \gamma_k^2\eta  \|u^k-u^{k+1}\|^2 -\gamma_k^2 \|K^*(u^k-u^{k+1})\|^2.\label{eq16}
\end{align}
Therefore, in the conditions of Theorem~\ref{th3}, let $u^\star$ be any solution of \eqref{eqpd}; that is, $u^\star\in\partial H(Kx^\star)$ and $0\in \partial R(x^\star)+\nabla F(x^\star)+K^* u^\star$. Then the constant $c_0$ is
\begin{equation*}
c_0=\frac{1-\gamma_0\mu_F\kappa}{\gamma_0^2}\Big(\|x^1-x^\star\|^2+\gamma_0^2\eta\|u^1-u^0\|^2-\gamma_0^2\|K^*(u^1-u^0)\|^2\Big)+ \eta\|u^1-u^\star\|^2.
\end{equation*}

The last statement in Theorem~\ref{th3} is obtained as follows. First, for every $k\geq 1$, 
$x_R^k = x^{k+1}- \gamma_k K^*(u^k-u^{k+1})$, so that $\|x_R^k-x^\star\|^2
\leq 2\|x^{k+1}-x^\star\|^2 + 2\|K\|^2\|\gamma_k(u^k-u^{k+1})\|^2$. Second, from \eqref{eq16}, $\|x^{k+1}-x^\star\|^2=O(1/k^2)$ and $(\eta -\|K\|^2)\|\gamma_k(u^k-u^{k+1})\|^2 \leq  \gamma_k^2 \langle (\eta I - KK^*) (u^k-u^{k+1}),u^k-u^{k+1}\rangle=O(1/k^2)$. So, assuming that $\eta>\|K\|^2$, $\|\gamma_k(u^k-u^{k+1})\|^2=O(1/k^2)$. Hence, $\|x_R^k-x^\star\|^2=O(1/k^2)$.

If $K=I$ and $\eta=1$, the PDDY algorithm reverts to the Davis--Yin algorithm, as given in \eqref{eqdybase}, but with $R$ and $H$ exchanged. In the conditions of Theorem~\ref{th3}, let $u^\star$ be any solution of \eqref{eqpd}; that is, $u^\star\in\partial H(x^\star)$ and $0\in \partial R(x^\star)+\nabla F(x^\star)+u^\star$. Then the constant $c_0$ is
\begin{equation*}
c_0=\frac{1-\gamma_0\mu_F\kappa}{\gamma_0^2}\|x^1-x^\star\|^2+ \|{\textstyle\frac{1}{\gamma_0}}(s^0-x^1)-u^\star\|^2.
\end{equation*}
This is the same value as in \eqref{eqc01}, corresponding to the Davis--Yin algorithm, viewed as the PD3O algorithm,  with $R$ and $H$ exchanged. Indeed, $u^\star$ is defined differently in both cases; that is, with the exchange, $u^\star\in \partial R(x^\star)$ in \eqref{eqc01}.

\subsection{$\boldsymbol{R=0}$: The Loris--Verhoeven Algorithm}

If $R=0$, the PD3O algorithm becomes:
Let $q^{0}\in\mathcal{X}$ and $u^{0}\in\mathcal{U}$. For $k=0,1,\ldots$ iterate:
\begin{equation}\left\lfloor\begin{array}{l}
x^{k+1} = \gamma_{k}(q^k-K^*u^{k})\\
q^{k+1}=\frac{1}{\gamma_{k+1}}x^{k+1} 
- \nabla F (x^{k+1})\\
u^{k+1} = \mathrm{prox}_{H^*/(\gamma_{k+1} \eta)} \big(u^k+\frac{1}{\eta}K(\frac{1}{\gamma_k}x^{k+1}+q^{k+1}-q^k )\big),
\end{array}\right.\label{eqlv100}\end{equation}
whereas the PDDY algorithm becomes: Let $x_R^{0}\in\mathcal{X}$, $u^{0}\in\mathcal{U}$. Set $p^{0}=K^*u^{0}$. For $k=0,1,\ldots$ iterate: \begin{equation*}\left\lfloor\begin{array}{l}
u^{k+1} = \mathrm{prox}_{H^*/(\gamma_{k} \eta)} \big(u^k +\frac{1}{\gamma_k \eta}Kx_R^k\big)\\
p^{k+1} = K^*u^{k+1}\\
x^{k+1}=x_R^k  - \gamma_k (p^{k+1}-p^k)\\
x^{k+1}_R = x^{k+1} -\gamma_{k+1}\nabla F (x^{k+1})- \gamma_{k+1} p^{k+1}.
\end{array}\right.\end{equation*}
Equivalently,
\begin{equation*}\left\lfloor\begin{array}{l}
u^{k+1} = \mathrm{prox}_{H^*/(\gamma_{k} \eta)} \big(u^k +\frac{1}{\gamma_k \eta}K(x^{k} -\gamma_{k}\nabla F (x^{k})- \gamma_{k} K^*u^k)\big)\\
x^{k+1}=x^{k} -\gamma_{k}\nabla F (x^{k})- \gamma_{k} K^*u^{k+1},
\end{array}\right.\end{equation*}
or:
\begin{equation*}\left\lfloor\begin{array}{l}
q^{k+1}=\frac{1}{\gamma_{k}}x^{k}- \nabla F (x^{k}) \\
u^{k+1} = \mathrm{prox}_{H^*/(\gamma_{k} \eta)} \big(u^k +\frac{1}{\gamma_k \eta}K(\gamma_{k}q^{k+1}- \gamma_{k} K^*u^k)\big)\\
x^{k+1}=\gamma_kq^{k+1}- \gamma_{k} K^*u^{k+1},
\end{array}\right.\end{equation*}
which is equivalent to \eqref{eqlv100}. So, when $R=0$, both the PD3O and PPDY revert to an algorithm which, for $\gamma_k\equiv \gamma$, is the Loris--Verhoeven algorithm \citep{lor11,com14,con19}.\medskip

Let $u^\star$ be any solution of \eqref{eqpd}; that is, $u^\star\in\partial H(Kx^\star)$ and $0\in \nabla F(x^\star)+K^* u^\star$. In the conditions of Theorem~\ref{th2}, $c_0$ is:
\begin{equation*}
c_0=\frac{1-\gamma_0\mu_F\kappa}{\gamma_0^2}\|x^1-x^\star\|^2+\|
q^0-{\textstyle\frac{1}{\gamma_0}}x^1-K^*(u^0-u^\star)
+\nabla F(x^\star)\|^2+\eta \|u^0-u^\star\|^2-\|K^*(u^0-u^\star)\|^2.
\end{equation*}
On the other hand, in Theorem~\ref{th3}, 
\begin{equation*}
c_0=\frac{1-\gamma_0\mu_F\kappa}{\gamma_0^2}\Big(\|x^1-x^\star\|^2+\gamma_0^2\eta\|u^1-u^0\|^2-\gamma_0^2\|K^*(u^1-u^0)\|^2\Big)+ \eta\|u^1-u^\star\|^2.
\end{equation*}
It is not clear how these two values compare to each other. They are both valid, in any case.

\subsection{$\boldsymbol{F=0}$: The Chambolle--Pock and Douglas--Rachford Algorithms}

If $F=0$, the PD3O algorithms reverts to:
Let $x^{0}\in\mathcal{X}$ and $u^{0}\in\mathcal{U}$. For $k=0,1,\ldots$ iterate:
\begin{equation*}\left\lfloor\begin{array}{l}
x^{k+1} = \mathrm{prox}_{\gamma_k R}\big(x^k-\gamma_{k}K^*u^{k} \big)\\
u^{k+1} = \mathrm{prox}_{H^*/(\gamma_{k+1} \eta)} \Big(u^k+\frac{1}{\eta}K\big((\frac{1}{\gamma_{k+1}}+\frac{1}{\gamma_k})x^{k+1}-\frac{1}{\gamma_{k}}x^{k} \big)\Big).
\end{array}\right.\end{equation*}

For $\gamma_k\equiv \gamma$, this is the form I \citep{con19} of the Chambolle--Pock algorithm \citep{cha11a}.

In the conditions of Theorem~\ref{th2}, let $u^\star$ be any solution of \eqref{eqpd}; that is, $u^\star\in\partial H(Kx^\star)$ and $0\in \partial R(x^\star)+K^* u^\star$. Then the constant $c_0$ is
\begin{equation*}
c_0=\frac{1}{\gamma_0^2}\|x^1-x^\star\|^2+\|
{\textstyle\frac{1}{\gamma_0}}(x^0-x^1)-K^*(u^0-u^\star)\|^2+\eta \|u^0-u^\star\|^2-\|K^*(u^0-u^\star)\|^2.
\end{equation*}

On the other hand, if $F=0$,  the PDDY algorithm reverts to: Let $x_R^{0}\in\mathcal{X}$, $u^{0}\in\mathcal{U}$. Set $p^{0}=K^*u^{0}$. For $k=0,1,\ldots$ iterate:
\begin{equation*}\left\lfloor\begin{array}{l}
u^{k+1} = \mathrm{prox}_{H^*/(\gamma_{k} \eta)} \big(u^k +\frac{1}{\gamma_k \eta}Kx_R^k\big)\\
p^{k+1} = K^*u^{k+1}\\
x^{k+1}=x_R^k  - \gamma_k (p^{k+1}-p^k)\\
x^{k+1}_R = \mathrm{prox}_{\gamma_{k+1} R}\big(x^{k+1}- \gamma_{k+1} p^{k+1}\big),
\end{array}\right.\end{equation*}
which can be simplified as: Let $x_R^{0}\in\mathcal{X}$, $u^{0}\in\mathcal{U}$. For $k=0,1,\ldots$ iterate:
\begin{equation*}\left\lfloor\begin{array}{l}
u^{k+1} = \mathrm{prox}_{H^*/(\gamma_{k} \eta)} \big(u^k +\frac{1}{\gamma_k \eta}Kx_R^k\big)\\
x^{k+1}_R = \mathrm{prox}_{\gamma_{k+1} R}\Big(x_R^k-  K^*\big((\gamma_k+\gamma_{k+1}) u^{k+1}-\gamma_k u^k\big)\Big),
\end{array}\right.\end{equation*}
knowing that we can retrieve the variable $x^k$ as $x^{k+1}=x_R^k  - \gamma_k K^*(u^{k+1}-u^k)$.

For $\gamma_k\equiv \gamma$, this is the form II \citep{con19} of the Chambolle--Pock algorithm \citep{cha11a}.

Note that with constant stepsizes, the Chambolle--Pock form II can be viewed as the form I applied to the dual problem. This interpretation does not hold with varying stepsizes as in Theorem~\ref{th2}: the stepsize playing the role of $\gamma_k$ would be $1/(\gamma_k\eta)$, which tends to $+\infty$ instead of 0, so that the theorem does not apply.

Note, also, that Theorem~\ref{th3} does not apply, since $F=0$ is not strongly convex. Finally, if the accelerated Chambolle--Pock algorithm form I is applied to the dual problem, our results do not guarantee convergence of the primal variable $x^k$ to a solution. So, we cannot derive an accelerated Chambolle--Pock algorithm form II.
\medskip

If $K=I$, $\mathcal{U}=\mathcal{X}$ and $\eta=1$, the Chambolle-Pock algorithm form I becomes the Douglas--Rachford algorithm: Let $x^{0}\in\mathcal{X}$ and $u^{0}\in\mathcal{X}$. For $k=0,1,\ldots$ iterate:
 \begin{equation*}\left\lfloor\begin{array}{l}
x^{k+1} = \mathrm{prox}_{\gamma_k R}\big(x^k-\gamma_{k}u^{k} \big)\\
u^{k+1} = \mathrm{prox}_{H^*/\gamma_{k+1} } \big(u^k+(\frac{1}{\gamma_{k+1}}+\frac{1}{\gamma_k})x^{k+1}-\frac{1}{\gamma_{k}}x^{k} \big).
\end{array}\right.\end{equation*}
We can rewrite the algorithm using only the meta-variable $s^k=x^k-\gamma_{k}u^{k}$: Let $s^{0}\in\mathcal{X}$. For $k=0,1,\ldots$ iterate:
 \begin{equation*}\left\lfloor\begin{array}{l}
x^{k+1} = \mathrm{prox}_{\gamma_k R}(s^k )\\
u^{k+1} = \mathrm{prox}_{H^*/\gamma_{k+1} } \big((\frac{1}{\gamma_{k+1}}+\frac{1}{\gamma_k})x^{k+1}-\frac{1}{\gamma_{k}}s^{k} \big)\\
s^{k+1}=x^{k+1}-\gamma_{k+1}u^{k+1}.
\end{array}\right.\end{equation*}
Using the Moreau identity, we obtain: Let $s^{0}\in\mathcal{X}$. For $k=0,1,\ldots$ iterate:
 \begin{equation}\left\lfloor\begin{array}{l}
x^{k+1} = \mathrm{prox}_{\gamma_k R}(s^k )\\
x_H^{k+1}= \mathrm{prox}_{\gamma_{k+1} H} \big((1+\frac{\gamma_{k+1}}{\gamma_k})x^{k+1}-\frac{\gamma_{k+1}}{\gamma_{k}}s^{k} \big)\\
s^{k+1}=x_H^{k+1}+\frac{\gamma_{k+1}}{\gamma_k}(s^{k}-x^{k+1}),
\end{array}\right.\label{eqdrapp1}\end{equation}
and for $\gamma_k\equiv \gamma$, we recognize the classical form of the Douglas--Rachford algorithm \citep{com10}.

In the conditions of Theorem~\ref{th2}, let $u^\star$ be any solution of \eqref{eqpd}; that is, $u^\star\in\partial H(x^\star)$ and $0\in \partial R(x^\star)+u^\star$. Then the constant $c_0$ is
\begin{equation*}
c_0=\frac{1}{\gamma_0^2}\|x^1-x^\star\|^2+\|
{\textstyle\frac{1}{\gamma_0}}(s^0-x^1)+u^\star\|^2.
\end{equation*}

On the other hand, if $K=I$, $\mathcal{U}=\mathcal{X}$ and $\eta=1$, the Chambolle-Pock algorithm form II becomes: 
Let $x_R^{0}\in\mathcal{X}$, $u^{0}\in\mathcal{U}$. For $k=0,1,\ldots$ iterate:
\begin{equation*}\left\lfloor\begin{array}{l}
u^{k+1} = \mathrm{prox}_{H^*/\gamma_{k}} \big(u^k +\frac{1}{\gamma_k }x_R^k\big)\\
 x^{k+1}=x_R^k  - \gamma_k (u^{k+1}-u^k)\\
x^{k+1}_R = \mathrm{prox}_{\gamma_{k+1} R}\big(x^{k+1}- \gamma_{k+1} u^{k+1}\big).
\end{array}\right.\end{equation*}
Using the Moreau identity, we obtain: 
Let $x_R^{0}\in\mathcal{X}$, $u^{0}\in\mathcal{U}$. For $k=0,1,\ldots$ iterate:
\begin{equation*}\left\lfloor\begin{array}{l}
x^{k+1}= \mathrm{prox}_{\gamma_{k}H} (x_R^k  + \gamma_k u^k)\\
u^{k+1} =u^k+(x_R^k   -x^{k+1})/\gamma_k\\
x^{k+1}_R = \mathrm{prox}_{\gamma_{k+1} R}\big(x^{k+1}- \gamma_{k+1} u^{k+1}\big).
\end{array}\right.\end{equation*}
Introducing the meta-variable $s^k=x_R^k+\gamma_{k}u^{k}$, we obtain: 
Let $s^{0}\in\mathcal{X}$. For $k=0,1,\ldots$ iterate:
\begin{equation*}\left\lfloor\begin{array}{l}
x^{k+1}= \mathrm{prox}_{\gamma_{k}H} (s^k)\\
x^{k+1}_R = \mathrm{prox}_{\gamma_{k+1} R}\big((1+\frac{\gamma_{k+1}}{\gamma_k})x^{k+1} - \frac{\gamma_{k+1}}{\gamma_k} s^k  \big)\\
s^{k+1}=x_R^{k+1}+\frac{\gamma_{k+1}}{\gamma_k}(s^k   -x^{k+1}).
\end{array}\right.\end{equation*}
Thus, we recover exactly  the Douglas--Rachford algorithm \eqref{eqdrapp1}, with $R$ and $H$ exchanged.

\section{Derivation of the Distributed Algorithms}\label{secA6}

\subsection{The Distributed PD3O Algorithm and its Particular Cases}

Let us adopt the notations of Section~\ref{secdi} and precise the different operators. The gradient of $\widehat{F}$ in $\widehat{\mathcal{X}}$ is
\begin{equation*}
\nabla\widehat{F}(\hat{x})=\big({\textstyle\frac{1}{M\omega_1}}\nabla F_1(x_1),\ldots,{\textstyle\frac{1}{M\omega_M}}\nabla F_M(x_M)\big),\quad \forall \hat{x}\in\widehat{\mathcal{X}}.
\end{equation*}
We define the linear subspace  
$\mathcal{S}=\{\hat{x}\in\widehat{\mathcal{X}}\ :\ x_1=\cdots=x_M\}$. $\widehat{F}$ is $L_{\widehat{F}}$-smooth, with 
$L_{\widehat{F}}=\max_m \frac{L_{F_m}}{M\omega_m}$. But since $\nabla \widehat{F}$ is applied to an element of $\mathcal{S}$ in the algorithms, we can weaken the condition on $L_{\widehat{F}}>0$ to be: 
 for every $\hat{x}=(x)_{m=1}^M \in \mathcal{S}$ and $\hat{x}'=(x')_{m=1}^M \in \mathcal{S}$,
\begin{align*}
\|\nabla\widehat{F}(\hat{x})-\nabla\widehat{F}(\hat{x}')\|^2_{\widehat{\mathcal{X}}}&=
\sum_{m=1}^M \omega_m \big\|{\textstyle\frac{1}{M\omega_m}}\nabla F_m(x)-{\textstyle\frac{1}{M\omega_m}}\nabla F_m(x')\big\|^2\\
&\leq L_{\widehat{F}}^2 \|\hat{x}-\hat{x}'\|^2_{\widehat{\mathcal{X}}}=L_{\widehat{F}}^2\|x-x'\|^2.
\end{align*}
That is, $L_{\widehat{F}}$ is 
such that, for every $(x,x')\in \mathcal{X}^2$,
\begin{equation}
\frac{1}{M^2}\sum_{m=1}^M \frac{1}{\omega_m}\|\nabla F_m(x)-\nabla F_m(x')\|^2\leq L_{\widehat{F}}^2 \|x-x'\|^2.\label{eqbeta}
\end{equation}
Notably, \begin{equation*}
L_{\widehat{F}}^2=\frac{1}{M^2}\sum_{m=1}^M \frac{L_{F_m}^2}{\omega_m}\end{equation*}
 satisfies the condition.

The adjoint operator of $\widehat{K}$ is
\begin{equation*}
\widehat{K}^*:\hat{u}\in\widehat{\mathcal{U}}\mapsto \big(K^*_1 u_1,\ldots,K^*_M u_M\big)\in\widehat{\mathcal{X}}.
\end{equation*}
Thus,
\begin{equation}
\|\widehat{K}\|^2=\|\widehat{K}^*\widehat{K}\|=\max_m \|K_m\|^2.\label{eqnorm1}
\end{equation}
But if $F_1=\cdots=F_M$, we can restrict the norm to $\mathcal{S}$ and 
\begin{align}
\|\widehat{K}\|^2&=\sup_{\hat{x}
\in \mathcal{S}}\, \langle \hat{x},\widehat{K}^*\widehat{K}\hat{x} \rangle_{\widehat{\mathcal{X}}}/\|\hat{x}\|^2_{\widehat{\mathcal{X}}}\notag\\
&= \sup_{x\in \mathcal{X}}\, \langle x, \sum_{m=1}^M \omega_m  K_m^*K_m x \rangle/\|x\|^2\notag\\
&= \big\|\sum_{m=1}^M \omega_m  K_m^*K_m\big\|,\label{eqnorm2}
\end{align}
which is $\leq \sum_{m=1}^M \omega_m \|K_m\|^2$.

For any $\zeta>0$, we have $\mathrm{prox}_{\zeta\widehat{R}}:\hat{x}\mapsto (x',\ldots,x')$, where $x'=\mathrm{prox}_{\zeta R}\big(\sum_{m=1}^M \omega_m x_m\big)$ and $\mathrm{prox}_{\zeta\widehat{H}}:\hat{u}\mapsto\big(\mathrm{prox}_{\zeta H_1/(M\omega_1)}(u_1),\ldots,\mathrm{prox}_{\zeta H_M/(M\omega_M)}(u_M)\big)$. We also have
 $\partial \widehat{H}:\hat{u}\mapsto \frac{1}{M\omega_1}\partial H_1(u_1) \times \cdots \times \frac{1}{M\omega_M}\partial H_M(u_M)$, 
$\widehat{H}^*:\hat{u}\mapsto \frac{1}{M}\sum_{m=1}^M H_m^*(M\omega_m u_m)$, and 
$\mathrm{prox}_{\zeta\widehat{H}^*}:\hat{u}\mapsto\big(\frac{1}{M\omega_1}\mathrm{prox}_{\zeta M \omega_1 H_1^*}(M\omega_1 u_1),$ $\ldots,\frac{1}{M\omega_M}\mathrm{prox}_{\zeta M \omega_M H_M^*}(M\omega_M u_M)\big)$.

By doing all these substitutions in the PD3O algorithm, we obtain the distributed PD3O algorithm, and all its particular cases, shown above. Theorem~\ref{th1} becomes Theorem~\ref{thb1} as follows. The objective function is $\Psi:x\in\mathcal{X}\mapsto R(x)+\frac{1}{M}\sum_{m=1}^M (F_m(x)+H_m(K_mx))$.

\begin{theorem}[convergence rate of the Distributed PD3O Algorithm]\label{thb1} In the Distributed PD3O Algorithm, suppose that $\gamma_k\equiv \gamma \in (0,2/L_{\widehat{F}})$, where $\widehat{F}$ satisfies \eqref{eqbeta}; if $F_m \equiv 0$, we can choose any $\gamma>0$. Also, suppose that $\eta\geq \|\widehat{K}\|^2$, where $\|\widehat{K}\|^2$ is defined in \eqref{eqnorm1} or \eqref{eqnorm2}. 
 Then $x^k$ converges to some solution $x^\star$ of \eqref{eqpbdi}. Also, $u_m^k$ converges to some element $u_m^\star\in\mathcal{U}_m$, for every $m=1,\ldots,M$. In addition, suppose that every $H_m$ is continuous on an open ball centered at $K_m x^\star$.  
 Then the following hold:
\begin{equation*}
\mathrm{(i)}\quad\Psi(x^k)-\Psi(x^\star) =o(1/\sqrt{k}).
\end{equation*}
Define the weighted ergodic iterate 
$\bar{x}^{k} = \frac{2}{k(k+1)}\sum_{i=1}^{k}i x^{i}$, for every $k\geq 1$. Then 
\begin{equation*}
\mathrm{(ii)}\quad\Psi(\bar{x}^k)-\Psi(x^\star) =O(1/k).
\end{equation*}
Furthermore, if every $H_m$ is $L_m$-smooth for some $L_m>0$, we have a faster decay for the best iterate so far:
\begin{equation*}
\mathrm{(iii)}\quad\min_{i=1,\ldots,k}\Psi(x^i)-\Psi(x^\star) =o(1/k).
\end{equation*}
\end{theorem}

The theorem applies to the particular cases of the Distributed PD3O Algorithm, like the distributed Loris--Verhoeven, Chambolle--Pock, Douglas--Rachford algorithms. We can note that the distributed forward--backward algorithm is monotonic, so Theorem~\ref{thb1} $\mathrm{(iii)}$ (with $H_m\equiv 0)$ yields $\Psi(x^k)-\Psi(x^\star) =o(1/k)$ for this algorithm. \bigskip

We now give accelerated convergence results using varying stepsizes, in presence of strong convexity. For this, we have to define the strong convexity constants $\mu_{\widehat{F}}$ and $\mu_{\widehat{R}}$. Like for the smoothness constant, we can restrict their definition to $\mathcal{S}$. So, $\mu_{\widehat{F}}$ becomes the strong convexity constant of the average function $\frac{1}{M} \sum_{m=1}^M F_m$. That is, $\mu_{\widehat{F}}\geq 0$ is such that the function
\begin{equation*}
x\in\mathcal{X} \mapsto \frac{1}{M} \sum_{m=1}^M F_m(x) - \frac{\mu_{\widehat{F}}}{2} \|x\|^2
\end{equation*} 
is convex. It is much weaker to require $\mu_{\widehat{F}}> 0$ than to ask all $F_m$ to be strongly convex. 
Similarly, we have $\mu_{\widehat{R}}=\mu_R$, the strong convexity constant of $R$. Thus, since the Accelerated Distributed PD3O Algorithm can be viewed as the accelerated PD3O algorithm applied to the minimization of $\widehat{F}(\hat{x})+\widehat{R}(\hat{x})+\widehat{H}(\widehat{K}\hat{x})$, we have all the ingredients to invoke Theorem~\ref{th2}, which is transposed as:

\begin{theorem}[Accelerated Distributed PD3O Algorithm]\label{thb2} Suppose that $\mu_{\widehat{F}}+\mu_R>0$.  Let $x^\star$ be the unique solution to \eqref{eqpbdi}. Let $\kappa\in(0,1)$ and $\gamma_0\in (0,2(1-\kappa)/L_{\widehat{F}})$.  Set $\gamma_1=\gamma_0$ and
\begin{equation*}
\gamma_{k+1}=\frac{-\gamma_k^2\mu_{\widehat{F}}\kappa+\gamma_k\sqrt{(\gamma_k\mu_{\widehat{F}}\kappa)^2+1+2\gamma_k\mu_R}}{1+2\gamma_k\mu_R}
,\quad\mbox{for every }k\geq 1.
\end{equation*}
Suppose that $\eta\geq \|\widehat{K}\|^2$, where $\|\widehat{K}\|^2$ is defined in \eqref{eqnorm1} or \eqref{eqnorm2}. 
Then in the Distributed PD3O Algorithm, there exists $\hat{c}_0>0$ such that, for every $k\geq 1$,
\begin{equation*}
\|x^{k+1}-x^\star\|^2\leq \frac{\gamma_{k+1}^2}{1-\gamma_{k+1}\mu_{\widehat{F}}\kappa}\hat{c}_0
=O\big(1/k^2\big).
\end{equation*}
\end{theorem}

As for Theorem~\ref{th4}, its counterpart in the distributed setting is:

\begin{theorem}[linear convergence of the Distributed PD3O Algorithm]\label{thb3} Suppose that $\mu_{\widehat{F}}+\mu_R>0$ and that every $H_m$ is $L_m$-smooth, for some $L_m>0$. Let $x^\star$ be the unique solution to \eqref{eqpbdi}. We suppose that $\gamma_k\equiv \gamma\in (0,2/L_{\widehat{F}})$ and $\eta\geq \|\widehat{K}\|^2$, where $\|\widehat{K}\|^2$ is defined in \eqref{eqnorm1} or \eqref{eqnorm2}. 
Then the Distributed PD3O Algorithm converges linearly: there exists $\rho\in (0,1]$ and $\hat{c}_0>0$ such that, for every $k\in\mathbb{N}$,
\begin{equation*}
\|x^{k+1}-x^\star\|^2\leq (1-\rho)^k \hat{c}_0.
\end{equation*}
\end{theorem}

We can remark that the Distributed Davis--Yin algorithm (with $\omega_m=1/M$ and $\gamma_k\equiv \gamma$) has been proposed in an unpublished paper by Ryu and Yin~\citep{ryu17}, where it is named Proximal-Proximal-Gradient Method. Their results are similar to ours in Theorems~\ref{thb1} and \ref{thb3} for this algorithm, but their condition $\gamma < 3/(2L)$, with $L=\max_m  L_{F_m}$, is worse than ours. Also, our accelerated version with varying stepsizes in Theorem~\ref{thb2} is new.

\subsection{The Distributed PDDY Algorithm}

The Distributed PDDY Algorithm, shown above, is derived the same way as the Distributed PD3O Algorithm. However, the smoothness constant cannot be defined only on $\mathcal{S}$, so that we have 
\begin{equation*}
L_{\widehat{F}}=\max_{m=1,\ldots,M} \frac{L_{F_m}}{M\omega_m}
\end{equation*}
and
\begin{equation*}
\mu_{\widehat{F}}=\min_{m=1,\ldots,M} \frac{\mu_{F_m}}{M\omega_m}.
\end{equation*}
Moreover, 
\begin{equation}
\|\widehat{K}\|^2=\max_{m=1,\ldots,M}  \|K_m\|^2,\label{eqkp1}
\end{equation}
except if $F_m \equiv 0$, in which case the Distributed PDDY Algorithm becomes the Distributed Chambolle--Pock Algorithm Form II, for which we can set  
\begin{equation}
\|\widehat{K}\|^2=\left\|\sum_{m=1}^M \omega_m  K_m^*K_m\right\|.\label{eqkp2}
\end{equation}

We can note that when $K_m\equiv I$, the Distributed PDDY Algorithm reverts to a form of distributed Davis--Yin algorithm, which is different from the Distributed Davis--Yin Algorithm obtained from the PD3O algorithm, shown above. Similarly, when $R=0$, we obtain a different algorithm than the Distributed Loris--Verhoeven Algorithm shown above. 
When $F_m\equiv 0$, the Distributed PDDY Algorithm reverts to the Distributed Chambolle--Pock Algorithm Form II, which is still different from the Distributed Douglas--Rachford Algorithm when $K_m\equiv I$.

The counterpart of Theorem~\ref{theoaj} is:
\begin{theorem}[convergence  of the Distributed PDDY Algorithm]\label{theoajdi}
In the Distributed PDDY Algorithm, suppose that $\gamma_k\equiv \gamma \in (0,2/L_F)$ and $\eta\geq \|\widehat{K}\|^2$, where $\|\widehat{K}\|^2$ is defined in \eqref{eqkp1} or \eqref{eqkp2}. Then all $x_m^k$ as well as   $x_R^k$ converge to the same solution $x^\star$  of \eqref{eqpbdi}, and every $u_m^k$ converges to some element $u_m^\star$.
\end{theorem}

The counterpart of Theorem~\ref{th3} is:
\begin{theorem}[Accelerated Distributed PDDY Algorithm]\label{thb4}Suppose that $\mu_{\widehat{F}}>0$. Let $x^\star$ be the unique solution to \eqref{eqpbdi}. Let $\kappa\in(0,1)$ and $\gamma_0\in (0,2(1-\kappa)/L_{\widehat{F}})$.  Set $\gamma_1=\gamma_0$ and
\begin{equation*}
\gamma_{k+1}=-\gamma_k^2\mu_{\widehat{F}}\kappa+\gamma_k\sqrt{(\gamma_k\mu_{\widehat{F}}\kappa)^2+1}
,\quad\mbox{for every }k\geq 1.
\end{equation*}
Suppose that $\eta\geq \|\widehat{K}\|^2$, where $\|\widehat{K}\|^2$ is defined in \eqref{eqkp1} or \eqref{eqkp2}. 
Then in the Distributed PDDY Algorithm, there exists $\hat{c}_0>0$ such that, for every $k\geq 1$,
\begin{equation*}
\sum_{m=1}^M \omega_m \|x_m^{k+1}-x^\star\|^2\leq \frac{\gamma_{k+1}^2}{1-\gamma_{k+1}\mu_F\kappa}c_0
=O\big(1/k^2\big).
\end{equation*}
Consequently, for every $m=1,\ldots,M$,
\begin{equation*}
 \|x_m^{k}-x^\star\|^2=O\big(1/k^2\big).
\end{equation*}
Moreover, if $\eta> \|\widehat{K}\|^2$, $\|x_R^{k}-x^\star\|^2 = O(1/k^2)$ as well.
\end{theorem}

The counterpart of Theorem~\ref{th4} is:
\begin{theorem}[linear convergence of the Distributed PDDY Algorithm]\label{thb5}Suppose that $\mu_{\widehat{F}}+\mu_R>0$ and that every $H_m$ is $L_m$-smooth, for some $L_m>0$. Let $x^\star$ be the unique solution to \eqref{eqpbdi}. Suppose that $\gamma_k\equiv \gamma \in (0,2/L_{\widehat{F}})$ and 
 $\eta\geq \|\widehat{K}\|^2$, where $\|\widehat{K}\|^2$ is defined in \eqref{eqkp1} or \eqref{eqkp2}. 
Then the Distributed PDDY Algorithm converges linearly: there exists $\rho\in (0,1]$ and $\hat{c}_0>0$ such that, for every $k\in\mathbb{N}$,
\begin{equation*}
\|x_R^{k+1}-x^\star\|^2\leq (1-\rho)^k \hat{c}_0.
\end{equation*}
\end{theorem}

\subsection{The Distributed Condat--V\~u Algorithm}

We can apply our product-space technique to other algorithms; in particular, we can derive distributed versions, shown below, of the Condat--V\~u algorithm \citep{con13,vu13, con19}, which is a well known algorithm for the problem \eqref{eqpb}.

The smoothness constant $L_{\widehat{F}}^2$ is the same as for the Distributed PD3O Algorithm; we can set $L_{\widehat{F}}^2=\frac{1}{M^2}\sum_{m=1}^M L_{F_m}^2/\omega_m$. 

Moreover,  the norm of $\widehat{K}$ is smaller for the Condat--V\~u algorithm: we have $\|\widehat{K}\|^2=\|\sum_{m=1}^M \omega_m  K_m^*K_m\|$, whatever the functions $F_m$. This is because the gradient descent step is completely decoupled from the dual variables in the Condat--V\~u algorithm.

\begin{figure*}[t!]
\begin{minipage}{.48\textwidth}
\begin{algorithm}[H]
		\caption*{\textbf{Distributed Condat--V\~u Alg.\ Form I}}
		\begin{algorithmic}
			\STATE \textbf{input:} $\gamma>0$, $\sigma>0$, $(\omega_m)_{m=1}^M$
			\STATE \ \ \ \ $x_0\in\mathcal{X}$, $(u_m^0)_{m=1}^M\in\widehat{\mathcal{U}}$%
			\STATE \textbf{initialize:} $a_m^{0}\coloneqq K_m^*u_m^{0}+\nabla F_m(x^{0})$, $\forall m$%
			\FOR{$k=0, 1, \ldots$}
			\STATE at master, \textbf{do}
			\STATE \ \ \ \ $x^{k+1} \coloneqq \mathrm{prox}_{\gamma R}\big(x^k-\frac{\gamma}{M}\sum_{m=1}^M a_m^{k} \big)$
			\STATE \ \ \ \ broadcast $x^{k+1}$ to all nodes
			\STATE at all nodes, for $m=1,\ldots,M$, \textbf{do}
			\STATE \ \ \ \ $u_m^{k+1} \coloneqq \mathrm{prox}_{M\omega_m \sigma H_m^*} \big(u_m^k$
			\STATE \ \ \ \ $\ \ \ \ {}+M\omega_m\sigma K_m(2 x^{k+1}-x^k )\big)$%
			\STATE \ \ \ \ $a_m^{k+1}\coloneqq  K_m^*u_m^{k+1}+\nabla F_m(x^{k+1})$ 
			\STATE \ \ \ \ transmit $a_m^{k+1}$  to master
			\ENDFOR
		\end{algorithmic}
	\end{algorithm}
	\end{minipage}
	\ \ \ \ \ \ \begin{minipage}{.48\textwidth}
\begin{algorithm}[H]
		\caption*{\textbf{Distributed Condat--V\~u Alg.\ Form II}}
		\begin{algorithmic}
			\STATE \textbf{input:} $\gamma>0$, $\sigma>0$, $(\omega_m)_{m=1}^M$
			\STATE \ \ \ \ $x_0\in\mathcal{X}$, $(u_m^0)_{m=1}^M\in\widehat{\mathcal{U}}$%
			\FOR{$k=0, 1, \ldots$}
			\STATE at all nodes, for $m=1,\ldots,M$, \textbf{do}
			\STATE \ \ \ \ $u_m^{k+1} \coloneqq \mathrm{prox}_{M\omega_m \sigma H_m^*} \big(u_m^k$
			\STATE \ \ \ \ $\ \ \ \ {}+M\omega_m\sigma K_m x^{k}\big)$%
			\STATE \ \ \ \ $a_m^{k}\coloneqq  K_m^*(2u_m^{k+1}-u_m^k)+\nabla F_m(x^{k})$%
			\STATE \ \ \ \ transmit $a_m^{k}$  to master
			\STATE at master, \textbf{do}
			\STATE \ \ \ \ $x^{k+1} \coloneqq \mathrm{prox}_{\gamma R}\big(x^k-\frac{\gamma}{M}\sum_{m=1}^M a_m^{k} \big)$
			\STATE \ \ \ \ broadcast $x^{k+1}$ to all nodes
			\ENDFOR
		\end{algorithmic}
	\end{algorithm}
	\end{minipage}
	\end{figure*}

The price to pay is a stronger condition on the parameters for convergence: 
\begin{theorem}[convergence of the Distributed Condat--V\~u Algorithm]\label{thb6}Suppose that the parameters $\gamma>0$ and $\sigma>0$ are such that 
\begin{equation*}
\gamma\Big(\sigma\big\|\sum_{m=1}^M \omega_m K_m^*K_m\big\|+\frac{L_{\widehat{F}}}{2}\Big)<1.
\end{equation*} 
Then $x^{k}$ converges to a solution $x^\star$ of \eqref{eqpbdi}. Also, $u_m^k$ converges to some element $u_m^\star\in\mathcal{U}_m$, for every $m=1,\ldots,M$.
\end{theorem}

When $F_m\equiv 0$, the two forms of the Distributed Condat--V\~u Algorithm revert to the two forms of the Distributed Chambolle--Pock Algorithm, respectively. In that case, with constant stepsizes $\gamma_k\equiv \gamma$, the convergence condition is $\gamma\sigma\|\sum_{m=1}^M \omega_m K_m^*K_m\|\leq 1$,  which is the same as above with $\sigma=1/(\eta\gamma)$.

\section*{Author Contributions}

Grigory Malinovsky wrote the code and generated the results for the SVM experiment in Section~\ref{sechl}. Peter Richt\'arik contributed to the paper writing and to the project management.  Laurent Condat did all the rest.

\bibliography{IEEEabrv,biblio}

\end{document}